\documentclass[opre,nonblindrev]{or_style/informs4_arxiv}

\OneAndAHalfSpacedXI 

\usepackage{endnotes}
\let\footnote=\endnote

\usepackage{graphicx}
\usepackage{joc_style/eqndefns-left}
\usepackage{multirow}
\usepackage{hhline}

\usepackage[small, margin=1cm]{caption}

\usepackage{color}
\definecolor{strcolor}{rgb}{0.6, 0.2, 0.6}
\definecolor{commentcolor}{rgb}{0.3125, 0.5, 0.3125}
\definecolor{keycol}{rgb}{0, 0, 1}

\usepackage{bbm}
\usepackage{booktabs}
\usepackage{algorithmic}
\usepackage{subfigure}
\usepackage{caption}
\usepackage{hyperref}
\usepackage{placeins}
\usepackage{pdflscape}
\usepackage{makecell}
\newcommand{\cell}[2]{\setlength{\tabcolsep}{0pt}\begin{tabular}{#1}#2 \end{tabular}}

\usepackage{listings}
\usepackage{url}

\usepackage{ulem}
\newcommand{\oo}[1]{{\textcolor{purple}{#1}}}
\newcommand{\stoo}[1]{\cl{\sout{#1}}}

\newcommand{\xxoo}[1]{{\stoo{[xxx]}}}

\newcommand{\cl}[1]{{\textcolor{blue}{#1}}}

\newcommand {\bea}{\begin{eqnarray}}
	\newcommand {\eea}{\end{eqnarray}}

\newcommand{\mb}[1]{\mbox{\boldmath \ensuremath{#1}}}

\newtheorem{algorithm}{Algorithm}

\newcommand{\firststage}{Group Representation Problem }
\newcommand{\fsabbrvx}{GRP}
\newcommand{\fsabbrv}{GRP }

\newcommand{\secondstage}{Assignment Problem under Fixed Representation }
\newcommand{\ssabbrv}{APFR }
\newcommand{\ssabbrvx}{APFR}

\renewcommand{\qed}{ \hfill\hbox{$\blacksquare$ }}
\def\blot{\quad \mbox{$\vcenter{ \vbox{ \hrule height.4pt
				\hbox{\vrule width.4pt height.9ex \kern.9ex \vrule width.4pt}
				\hrule height.4pt}}$}}

\usepackage{natbib}
\bibpunct[, ]{(}{)}{,}{a}{}{,}%
\def\bibfont{\fontsize{8}{9.5}\selectfont}%
\TheoremsNumberedThrough     
\ECRepeatTheorems 

\EquationsNumberedThrough    

\gdef\AQ#1{}
\gdef\CQ#1{}
\begin{document}
	
\def\COPYRIGHTHOLDER{INFORMS}%
\def\COPYRIGHTYEAR{2017}%
\def\DOI{\fontsize{7.5}{9.5}\selectfont\sf\bfseries\noindent https://doi.org/10.1287/opre.2017.1714\CQ{Word count = 9740}}

	\RUNAUTHOR{Lawless and G\"unl\"uk} %

	\RUNTITLE{Fair Clustering with Minimum Representation Constraints}

\TITLE{Fair Clustering with Minimum Representation Constraints 
}


	\ARTICLEAUTHORS{

\AUTHOR{Connor Lawless\textsuperscript{a}, Oktay G\"unl\"uk\textsuperscript{b}}
\AFF{$^{a}$Department of Management Science and Engineering, Stanford University\\
$^{b}$H. Milton Stewart School of Industrial and Systems Engineering, Georgia Institute of
Technology }


}
	


\ABSTRACT{Clustering is a well-studied unsupervised learning task that aims to partition data points into a number of clusters. In many applications, these clusters correspond to real-world constructs (e.g., electoral districts, playlists, TV channels), where a group (e.g., social or demographic) benefits only if it reaches a minimum level of representation in the cluster (e.g., 50\% to elect their preferred candidate). 
In this paper, we study the $k$-means and $k$-medians clustering problems under the additional fairness constraint that each group must attain a minimum level of representation in at least a specified number of clusters. We formulate this problem as a mixed-integer (nonlinear) optimization problem and propose an alternating minimization algorithm, called \textit{MiniReL}, to solve it.
Although incorporating fairness constraints results in an NP-hard assignment problem within the MiniReL algorithm, we present several heuristic strategies that make the approach practical even for large datasets. Numerical results demonstrate that our method yields fair clusters without increasing clustering cost across standard benchmark datasets.}



\KEYWORDS{Clustering, Integer Programming, Lloyd's Algorithm, Fair Machine Learning}

	
	%
	
\maketitle
	
\section{Introduction}

Clustering is an unsupervised learning task that aims to partition data points into sets of similar data points called clusters \citep{xu2005survey}. Clustering is widely used across various domains including customer segmentation \citep{kansal2018customer}, grouping content in entertainment platforms \citep{daudpota2019video}, and identifying subgroups within a clinical study \citep{wang2020unsupervised}, among others. 
However the widespread application of clustering, and machine learning more broadly, to human-centric tasks has raised concerns about its disparate impact on minority groups and other vulnerable demographics. Motivated by a flurry of recent results highlighting bias in many automated decision making tasks such as facial recognition \citep{buolamwini2018gender} and criminal justice \citep{mehrabi2019survey}, researchers have begun focusing on mechanisms to ensure that machine learning algorithms are \textit{fair} to all those affected. One of the challenges of fairness in an unsupervised learning context, compared to the supervised setting, is the lack of ground truth labels. 
Consequently, instead of enforcing approximately equal error rates across groups, fair clustering generally aims to ensure that composition of the clusters or their centers  (for settings like $k$-means and $k$-medians clustering)  represent all groups equitably.

A common approach to fair clustering is to enforce that each cluster includes data points from each group in proportion to their presence in the overall dataset
(i.e., via balance \citep{chierichetti2017fair} or bounded representation \citep{ahmadian2019clustering}). This approach aims to balance the presence of each group in each cluster and therefore tries to spread each group uniformly across the clusters.
Notice that this approach might not be desirable in settings where a group only gains a significant benefit from the cluster when they reach a minimum level of representation in that cluster. Consider the problem of clustering a set of media (e.g., TV shows) into cohesive segments (e.g., channels). A natural fairness consideration in designing these segments is to ensure that there is sufficient representation for different demographic groups. In these settings the benefit of the representation is only realized when a large percentage of the segment is associated with a demographic group (i.e., so that viewers can consistently watch programming that speaks to them). This is even legislated in some countries, for example Canadian television channels are required to have at least 50\% Canadian programming \citep{crtc}. Note that, in this setting, spreading a minority demographic group across all clusters ensures that the demographic group will never have majority representation in any segment.
As another example, consider a simple voting system for a committee where the goal is to first cluster voters (e.g., employees, faculty) into different constituencies that can then elect a committee representative. Here, a proportionally fair clustering would assign a minority group that represents 30\% of the vote equally among each cluster. However, the minority group only gets a benefit (i.e., the ability to elect a candidate of their choice) if they have at least 50\% representation in the cluster. 

In this paper we introduce a new notion of fairness in clustering that addresses this issue. Specifically, we introduce \textit{minimum representation fairness} (MR-fairness) which requires each group to have a certain number of clusters where they cross a given minimum representation threshold (i.e., 50\% in the voting example).

Arguably the most popular algorithms for clustering is Lloyd's algorithm for $k$-means clustering \citep{lloyd1982least}, and the associated alternating minimization approach for $k$-medians \citep{park2009simple}. 
These iterative algorithms alternate between fixing cluster centers and assigning each point to the closest cluster center.
Both algorithms are guaranteed to return a locally optimal
solution (i.e., no perturbation of the cluster centers around the solution leads to a better clustering cost). 
However, using these algorithms can lead to clusters that violate MR-fairness. 

In this paper we introduce a modified version of Lloyd's algorithm for $k$-means and $k$-medians that ensures minimum representation fairness, henceforth referred to as MINImum REpresentation fair Lloyd's algorithm (MiniReL for short). The key modification behind our approach is to replace the original greedy assignment step with a new optimization problem that finds the minimum cost assignment while ensuring fairness. In contrast to the standard clustering setting, we show that finding a minimum cost clustering that respects MR-fairness is NP-Hard even when the cluster centers are already fixed. We show that this optimization problem can be solved via integer programming (IP) in practice and introduce a number of computational approaches to improve the run-time. We empirically show that our approach is able to construct fair clusters which have nearly the same clustering cost as those produced by Lloyd's algorithm.

\subsection{Related Work}
A recent flurry of work in fair clustering has given rise to a number of different notions of fairness. One broad line of research, started by the seminal work of \citet{chierichetti2017fair}, puts constraints on the \textit{proportion} of each cluster that comes from different groups. This can be in the form of balance \citep{chierichetti2017fair,bera2019,schmidt2018fair,bercea2018cost,backurs2019scalable,kleindessner2019guarantees,ahmadian2020fair,bohm2020fair,chhabra2020fair,liu2021stochastic,ziko2021variational,le2021fair} which ensures each group has relatively equal representation, or a group specific proportion such as the bounded representation criteria \citep{ahmadian2019clustering,bera2019,ahmadian2020fair,schmidt2018fair,esmaeili2020probabilistic,jia2020fair,huang2019coresets,bandyapadhyay2020coresets,harb2020kfc} or maximum fairness cost \citep{chhabra2020fair}. MR-fairness bares a resemblance to this line of work as it puts a constraint on the proportion of a group in a cluster, however instead of constraining a fixed proportion across all clusters it looks holistically across all clusters and ensures that threshold is met in a baseline number of clusters. Another line of work tries to minimize the \textit{worst case average clustering cost} (i.e., $k$-means cost) over all the groups, called social fairness \citep{ghadiri2021socially,abbasi2021fair,makarychev2021approximation,goyal2021tight}. 
{Most similar to our algorithmic approach is the Fair Lloyd algorithm introduced in \citep{ghadiri2021socially} which studies social fairness. They also present a modified version of Lloyd's algorithm that converges to a local optimum. As a consequence of the social fairness criterion their approach requires a modified center computation step that can be done in polynomial time. MR-fairness, however, requires a modified cluster assignment step that is NP-hard which we solve via integer programming.

Most similar to MR-fairness is diversity-aware fairness introduced in \citep{thejaswi2021diversity} and the related notion of fair summarization \citep{kleindessner2019fair,chiplunkar2020solve,jones2020fair}. Under these fairness criteria, each group must have at least a specified number of its members included among the cluster centers.
MR-fairness differs in that our criteria is not tied to the group membership of the cluster center selected but the proportion of each group in a given cluster. 
Our fairness criterion
is more relevant in settings where the the demographic composition of each cluster impacts the societal outcome (e.g., voting) as opposed to settings where the cluster center itself impacts the outcome (e.g., committee selection where centers correspond to the committee chair).

There is also a long line of research that looks at fairness in the context of gerrymandering \citep{kueng2019fair,gurnee2021fairmandering,benade2022political,levin2019automated,mehrotra1998optimization,ricca2013political}. While our fairness criterion shares some similarity with different fairness criteria in gerrymandering, the gerrymandering problem places different constraints on the construction of the clusters such as contiguity. Consequently, algorithmic approaches to tackle gerrymandering require more computationally intensive optimization procedures.

Existing work has leveraged exact optimization algorithms to find globally optimal solutions to constrained clustering problems using tools such as convex optimization \citep{chen2014clustering}, column generation \citep{babaki2014constrained}, and IP \citep{agoston2024mixed, piccialli2022exact, miyauchi2018exact}. While these methods can be highly accurate, they are computationally intensive and often struggle to scale beyond 1,000 data points.
Recent research has successfully used IP as a subroutine within heuristic clustering algorithms to enforce constraints such as must-link and cannot-link constraints \citep{baumann2024algorithm}, cardinality constraints \citep{mueller2010integer, tang2019size}, and proportional fairness \citep{tran2023mpfcc}. Our work builds on these approaches by not only using IP as a sub-routine but also incorporating optimization-informed approximation techniques to scale the approach to much larger datasets (over $90,000$ data points). 
\subsection{Main Contributions}

We summarize our main contributions as follows:
\begin{itemize}
\item We introduce a novel definition of fairness for clustering of practical importance called MR-fairness, which requires that a specified number of clusters should have at least $\alpha$ percent members from a given group.
\item We formulate the MR-fair $k$-means and $k$-medians clustering as a mixed-integer optimization problem and introduce a new algorithm MiniReL, based on Lloyd's algorithm for clustering, to find a local optimum.

\item We prove that MR-fairness is incompatible with existing fairness criteria in clustering, underscoring the need for new algorithms to address it.

\item We show that unlike other definitions of fairness, MR-fairness can not be approximated by simply adjusting unfair cluster centers.

\item We show that incorporating MR-fairness into Lloyd's algorithm leads to a NP-Hard sub-problem to assign data points to fixed cluster centers, which we call the Fair Minimum Representation Assignment (FMRA) Problem. 

\item We introduce a two-stage approach to solve the FMRA problem using a polynomial-time bi-criteria approximation algorithm based on a network flow formulation. We also introduce a fast heuristic for setting the first-stage variables.

\item We present numerical results to demonstrate that MiniReL is able to construct MR-fair clusters with only a modest increase in run-time and with little to no increase in clustering cost compared to the standard $k$-means or $k$ -medians clustering algorithm.  
\end{itemize}


The remainder of the paper is organized as follows. In Section \ref{sec:MRFS} we formally describe the the MR-fair clustering problem and highlight its differences from other notions of fairness. In Section \ref{sec:mio} we present a mixed integer optimization formulation for the MR-fair clustering problem and introduce MiniReL. In Section \ref{sec:scaling} we introduce computational approaches to help our algorithm scale to large datasets. Finally Section \ref{sec:exp} presents a numerical study of MiniReL compared to the standard Lloyd's algorithm.

\section{Minimum Representation Fair Clustering}\label{sec:MRFS}
The input to the standard clustering problem is a set of $n$ $m$-dimensional data points $\mathcal{X} = \{x^i \in \mathbbm{R}^m\}_{i=1}^n$. Note that assuming the data points to have real-valued features is not a restrictive assumption as categorical features can be converted to real-valued features through the one-hot encoding scheme. The goal of the clustering problem is to partition the data points into a set of $K$ clusters $\mathcal{C} = \{C_1, \dots, C_K \}$ in such a way that some measure of clustering cost is minimized.
Let ${\cal K} = \{1, \dots, K\}$ be the index set of clusters. 

In the fair clustering setting, each data point also has a (small) number of sensitive features such as gender and race associated with it. Each sensitive feature can take a finite number of possible values (e.g., male, female, or non-binary for gender). Note that a sensitive feature could also correspond to an \textit{intersection of characteristics} (i.e., to promote intersectional fairness \citep{gohar2023survey}), for example we may have a sensitive feature might correspond to race \textit{and} gender (e.g., Black men, White men, Black women, White women). The choice of whether to consider intersectional fairness or each feature independently should be chosen based on the application.
We denote the set of sensitive features with ${\cal F}$ and possible values of a feature $f \in {\cal F}$ with the set ${\cal G}_f$.
We use ${\cal G} = \cup_{f \in {\cal F}} {\cal G}_f$ to be the set of all possible values of all features, and call each $g\in{\cal G}$ a group. 
With this notation, each data point $x^i$ is associated with $|{\cal F}|$ many groups, one for each $f\in \mathcal{F}$ (i.e., gender, race).
Let $X_g$ be the set of data points associated with group $g\in \mathcal{G}$ and note that unlike other fair machine learning work, we do not assume that $\{X_g\}_{g\in{\cal G}}$ form a partition of $\mathcal{X}$ when $|{\cal F}|>1$. 
Instead, $\{X_g\}_{g\in{\cal G}_f}$ forms a partition of $\mathcal{X}$ for each sensitive feature  $f \in {\cal F}$. For ease of reference, we include a table of common notation used in this paper in Appendix \ref{app:notation}.

The key intuition behind MR-fair clustering is that individuals belonging to a group gain a benefit only when their representation within a cluster exceeds a minimum threshold. We denote this minimum representation threshold with $\alpha \in (0,1]$, and define $\alpha$-representation as follows:

\begin{definition}[$\alpha$-representation]\label{def:alpha}
   A group $g\in \mathcal{G}$  is said to be $\alpha$-represented in a cluster $C_k$ if 
    $$
    | C_k \cap X_g | \geq \alpha |C_k|
    $$
    \end{definition}
    
The parameter $\alpha$ should be selected based on the the specific application. For example, in many voting systems a group must constitute a majority to influence outcomes (i.e., $\alpha = 0.5$). Our framework also allows for $\alpha$ to be group-dependent (i.e., a different $\alpha_g$ for each group $g$), however throughout the paper we assume $\alpha$ to be the same for all groups.
For a given clustering $\mathcal{C}$, group $g\in{\cal G}$, and $\alpha \in (0,1]$, let $\Lambda(\mathcal{C}, X_g, \alpha)$ be the number of clusters where group $g$ has $\alpha$-representation. In MR-fairness, each group $g$ has a parameter $\beta_g$ that specifies a minimum number of clusters where that group should have $\alpha$-representation.

\begin{definition}[Minimum representation fairness]
A given clustering $\mathcal{C} = \{C_1, \dots, C_K\}$ is said to be ($\alpha,\boldsymbol{\beta}$)-minimum representation fair  if for every group $g \in {\mathcal{G}}$:
\begin{equation}\Lambda(\mathcal{C}, X_g, \alpha) \geq \beta_{g} \label{eq:fairness_constraint}\end{equation}
for a given $\boldsymbol{\beta} = \{\beta_g \in \mathbb{Z}^+\}_{g \in \mathcal{G}}$.
\end{definition}

The definition of MR-fairness is flexible enough that the choice of  $\boldsymbol{\beta}$ can (and should) be specialized to each application together with the choice of $\alpha$. In Section \ref{sec:select_params}, we discuss how to set these parameters to reflect two common fairness criteria from the fair classification literature. Also note that at most $\lfloor 1/\alpha\rfloor$ groups corresponding to a single feature $f$ can be $\alpha$-represented in a cluster.

We integrate MR-fairness into two popular clustering paradigms: $k$-means and $k$-medians.
In both cases, the goal is to find a clustering along with a set of cluster centers $\{c_k \in {\cal L}_k \}_{k=1}^{K}$ that minimizes the clustering cost. In the $k$-means setting the center of a cluster can be located anywhere in the feature space (i.e., ${\cal L}_k = \mathbbm{R}^m ~~\forall k$), whereas in $k$-medians it must be at one of the data points assigned to the cluster (i.e., ${\cal L}_k = C_k$). Note that in the standard $k$-medians setting it is not necessary to explicitly require $c_k\in C_k$ as any (local) optimal solution would  satisfy this condition. However, this is no longer the case in the MR-fairness setting and therefore we explicitly constrain each center to belong to its corresponding cluster. Note that the $k$-median problem described in this paper is also referred to as the $k$-medoids clustering problem in parts of the machine learning community \citep{park2009simple}. We denote the cost of assigning a data point $x^i$ to a cluster with center $c$ with $D(x^i, c)$ and the objective is to minimize the total clustering cost over all data points.
In the $k$-means setting the cost is equal to the squared distance between the data point and the center (i.e., $D(x^i, c) = \|x^i - c\|^2$), whereas in $k$-medians one can use any distance function. 
Combining the standard clustering problem with MR-fairness requirement leads to the following optimization problem:

\begin{definition}[Minimum representation fair $k$-means and $k$-medians problem]
For a given $\alpha \in (0,1]$ and $\boldsymbol{\beta} = \{\beta_g \in \mathbb{Z}^+\}_{g \in \mathcal{G}}$, the minimum representation fair $k$-means and $k$-medians problems are:
$$
\min_{\mathcal{C}, c_1,\ldots,c_K} \sum_{k \in \mathcal{K}}\sum_{x^i \in C_k} D(x^i, c_k)
\quad \textbf{s.t.} \quad 
\Lambda(\mathcal{C}, X_g, \alpha) \geq \beta_{g} ~~\forall g \in \mathcal{G},\quad
c_k \in {\cal L}_k ~~ \forall k \in {\cal K}
$$
where ${\cal L}_k = \mathbbm{R}^m ~~ \forall k$, $D(x^i, c) = \|x^i - c\|_2^2$~~for $k$-means, and ${\cal L}_k = C_k$, $D(x^i, c)$ is any given distance function for  $k$-medians problem.

\end{definition}
Note that one important difference between the  fair  and the standard versions of the $k$-means and $k$-medians clustering problem is that greedily assigning data points to their closest cluster center may no longer be feasible for the fair version (i.e., assigning some data points to farther cluster centers may be necessary to meet the fairness criteria). Thus the problem can no longer be viewed as an optimization problem over cluster centers.

\subsection{Comparison to existing fair clustering criteria} \label{sec:fair_metric_comp}

In this section, we compare MR-fairness to existing fairness criteria in the clustering literature and show that it is fundamentally different. Specifically, we demonstrate that MR-fairness is \textit{incompatible} with both proportional and social-fairness, meaning that these criteria of fairness cannot be simultaneously satisfied. This implies that existing algorithms for fair clustering could fail to achieve MR-fairness, which we also show empirically in Section \ref{sec:exp}.

We start by formalizing the concept of incompatibility between two fairness criteria in clustering. We say two fairness criteria are incompatible if there exists a clustering instance in which both criteria cannot be satisfied simultaneously. Since many fairness criteria involve parameters (e.g., $\alpha$, $\mathbf{\beta}$ in MR-fairness), compatibility may depend on the choice of these parameters. To address this, we adopt a strict notion of incompatibility and say that two criteria are incompatible if there exists no set of non-trivial parameters that allow both criteria to be met at the same time.

Let $f_{\mathbf{p}}: {\cal C} \rightarrow \{0,1\}$ be a fairness criterion, parametrized by ${\mathbf{p}}$, which maps a clustering ${\cal C}$ to a binary response that represents whether or not it satisfies the criterion. Let $\Psi({\cal X}, {\cal G}, f_{\mathbf{p}})$ be the set of all clusterings that satisfy fairness criterion $f_{\mathbf{p}}$ for a dataset ${\cal X}$ with groups ${\cal G}$. With this notation in mind, we define incompatibility between two fairness criteria $f$ and $g$ as follows:

\begin{definition}
    A parametrized fairness criterion $f$ is said to be \textit{incompatible} with another criterion $g$ if there exists a fair clustering instance ${\cal X}, ~{\cal G}$ and fairness parameter $\mathbf{p}_1$ for $f$ such that 
    $$\Psi({\cal X}, {\cal G},f_{\mathbf{p}_1}) ~\bigcap~ \Psi({\cal X}, {\cal G},g_{\mathbf{p}_2}) = \emptyset$$
    for all non-trivial $\mathbf{p_2}$ (i.e., parameters that are not satisfied by all clusterings of the dataset). 
\end{definition}

Note that the definition of incompatibility is \textit{directional} --- a fairness criterion $f$ can be incompatible with a fairness criterion $g$ while the reverse might not be true. 
We next highlight the difference between MR-fairness and other well-known fairness criteria.

\noindent
\subsubsection{Proportional Fairness:} Several fairness criteria in the clustering literature (e.g., balance, bounded representation, group fairness) put constraints on the proportion of data points in each cluster that belong to different groups:

\begin{definition}
    A clustering $C$ for a dataset ${\cal X}$, ${\cal G}$ is said to be \textit{proportionally fair} (also known as group fairness) if $\forall k \in {\cal K}, g \in {\cal G}$:
    $$
    p_{lg}|C_k| \leq |C_k ~\cap ~X_g| \leq p_{ug}|C_k|
    $$
    where $p_{lg}, p_{ug}$ define the lower and upper bounds for the proportion of individuals belonging to group $g$ in every cluster.
\end{definition}
Note that this definition applies to \textit{every} cluster and thus ensures a similar demographic breakdown across all clusters. In contrast, MR-fairness may require clusters to have completely different demographic breakdowns (e.g., include a cluster where a minority group has majority representation). The following proposition proves that these two fairness criteria are in fact incompatible and cannot be guaranteed to be satisfied together:
\begin{proposition}
MR-fairness is incompatible with proportional fairness.
\end{proposition}
\proof{Proof}
Consider a simple clustering instance with two groups each with two data points and $K=2$, and consider MR-fairness constraints with $\alpha > 0.5$ and $\beta_g = 1$ for both groups. To achieve MR-fairness at least one cluster has to be homogeneous (i.e., all data points must belong to the same group). Note that non-trivial parameters for proportional fairness must either have $p_{lg} > 0$ or $p_{ug} < 1$ (i.e., there is at least some constraint on the proportion of each group in each cluster). Since one cluster in any MR-fair solution is entirely constructed from one group, there is no non-trivial setting of the proportional fairness parameters that can be simultaneously satisfied.
\qed

\noindent
\subsubsection{Social Fairness} Social fairness aims to minimize the average clustering cost of the worst-case group. Similar to prior work, we consider a constrained version of social fairness that requires every group to have a cost within a multiplicative factor of the optimal socially fair solution.

\begin{definition}
For a dataset ${\cal X}$ with groups ${\cal G}$ and a given $K$, a clustering $\bar C$ with centers $\bar c$ is said to be \textit{socially fair} if:
    $$
    \max_{g \in {\cal G}} \frac{1}{|X_g|}\sum_{k \in {\cal K}} \sum_{x^i \in X_g \cap \bar C_k}  D(x^i, \bar c_k) 
    \leq \alpha_{SF} \min_{{\cal C}, c_1, \dots, c_K} \max_{g \in {\cal G}} \frac{1}{|X_g|}\sum_{k \in {\cal K}} \sum_{x^i \in X_g \cap C_k}  D(x^i, c_k) 
    $$
    for a given fairness parameter $\alpha_{SF} \geq 1$.
\end{definition}

\begin{proposition} \label{prop:social_fairness}
MR-fairness is incompatible with social fairness
\end{proposition}
\proof{Proof}
Consider a simple 1-dimensional clustering problem with $K=2$ and two groups, each with two data points. For both groups one data point is located at $0$ and the other is located at $\gamma$. Note that the optimal solution to the unconstrained clustering problem is to place one center at $0$ and one center at $\gamma$ which has cluster cost $0$ for both groups (i.e., the optimal socially fair solution has cost $0$). Now consider an MR-fair constraint that requires each group to have at least one $\alpha$-represented cluster (i.e., $\beta = 1$) for $\alpha > 0.5$. Note that any MR-fair clustering requires assigning at least one data point to a center which has non-zero cost (e.g., keep centers at $0$ and $\gamma$ and assign one data point at $0$ to the center at $\gamma$). This violates social fairness for any choice of $\alpha_{SF}$ or any additive violation of the constraint less than $\gamma$. 
\qed

\noindent
\subsubsection{Diversity-Aware Fairness:} Designed specifically for the $k$-median problem, diversity-aware fairness takes a holistic view of group representation across the entire clustering. However unlike MR-fairness, it considers only the demographic composition of the cluster centers. Let $\bar{\cal L} \subseteq {\cal X}$ denote the set of selected centers for the $k$-median problem. 

\begin{definition}
    A clustering $C$ for a dataset ${\cal X}$, ${\cal G}$ is said to be \textit{diversity aware} if $g \in {\cal G}$:
    $$
    k_{lg} \leq |\bar{{\cal L}} ~\cap ~X_g| \leq k_{ug}
    $$
    where parameters $k_{lg}$ and $k_{ug}$ define lower and upper bounds for the number of cluster centers that belong to group $g$.
\end{definition}
Unlike the previous fairness criteria, MR-fair clustering is \textit{compatible} with diversity-aware clustering. 
Specifically, for any MR-fair clustering problem, we can construct a compatible diversity-aware formulation by setting $k_{lg} = \beta_g$ and $k_{ug} = K$. Any clustering that satisfies MR-fair clustering can be adjusted to satisfy diversity-aware fairness by simply adjusting the cluster centers to meet the diversity-aware criteria. This is always possible as any cluster that has $\alpha$-representation for a group must have at least one data point from that group which can be selected as a center. However, the reverse is not true and a diversity-aware clustering can have an arbitrarily large violation of MR-fairness.


\begin{proposition}
A diversity-aware fair clustering can violate MR-fairness for all groups.
\end{proposition}
\proof{Proof}
Consider the same simple clustering instance and MR-fair clustering constraint as Proposition \ref{prop:social_fairness}. We now construct a diversity-aware fairness constraint by setting $k_{lg} = 1$ for both groups and $k_{ug} = 2$. Note that the optimal clustering without fairness constraints (i.e., centers at $0$ and $\gamma$) can be made diversity-aware by selecting one center from each group. However, this clustering violates MR-fairness for both groups. Note that this example can be extended to have more groups by replicating the construction with additional groups each consisting of two data points.
\qed

Note that these incompatibility results do not rule out the possibility of satisfying multiple fairness criteria simultaneously in specific cases. However, they imply that no algorithm can produce a clustering that simultaneously satisfies incompatible fairness criteria. In Section \ref{sec:exp} we show that these incompatibilities do arise in practice, underscoring the need for new algorithms to ensure MR-fairness.

\subsection{Setting Fairness Parameters} \label{sec:select_params}

MR-fairness parameters $\boldsymbol{\beta}$ and $\alpha$ should be tailored to the specific application.
In this paper we consider two natural choices for $\boldsymbol{\beta}$ that mirror fairness definitions from the fair classification literature \citep{darlington1971another, hardt2016equality}. The first one is called \textit{cluster statistical parity} and it sets $\beta_g$ to be equal for all groups. Similar to its analog in fair classification, \textit{statistical parity} 
\citep{darlington1971another}, this choice of $\boldsymbol\beta$ is suitable for settings where we expect each group to have equal representation (e.g., male and female loan applicants).
\begin{definition}[Cluster Statistical Parity]
A clustering is said to satisfy \textit{cluster statistical parity} if it is a minimum representation fair clustering for a given $\alpha$ and the following $\beta$:
$$
\beta_g = \Big\lfloor \frac{1}{|\mathcal{G}_f|}\big\lfloor \alpha^{-1} \big \rfloor K \Big \rfloor \quad \forall g \in \mathcal{G}, f \in {\cal F}
$$
\end{definition}
Motivated by the notion of \textit{equality of opportunity} \citep{hardt2016equality} in supervised learning, the second choice is called \textit{cluster equality of opportunity} and it  sets $\beta_g$ to be proportional to the size of the group. 
\begin{definition}[Cluster Equality of Opportunity]
A clustering is said to satisfy \textit{cluster equality of opportunity} if it is a minimum representation fair clustering for a given $\alpha$ and the following $\boldsymbol\beta$:
$$
\beta_g = \Big \lfloor \frac{|X_g|}{n} \big \lfloor \alpha^{-1} \big\rfloor K \Big \rfloor \quad \forall g \in \mathcal{G}$$
\end{definition}
This setting of $\mathbf{\beta}$ is suitable for applications where the size of groups is imbalanced. For instance, in a congressional districting example if a minority group only represents 10\% of the population, it should not have majority voting power in 50\% of the districts.
We emphasize that although we use these fairness criteria in our computational experiments, it is important to select $\boldsymbol\beta$ and $\alpha$ that align with the intended application. For instance, in congressional districting $\boldsymbol\beta$ should 
reflect legal requirements such as those set by the Voting Rights Act, which may mandate  
majority-minority districts for protected groups (e.g., a majority-Black voting district). We also note that the MR-fair clustering problem may be infeasible for some choices of the parameters $\boldsymbol\beta$ and $\alpha$. In Section \ref{sec:feasibility} we introduce a quick method to check the feasibility of a given parameter setting.

\section{Mixed Integer Optimization Framework} \label{sec:mio}
We next formulate the MR-fair clustering problem as a mixed-integer program with a non-linear objective. We use binary variable $z_{ik}$ to denote if data point $x^i$ is assigned to cluster $k$, and variable $c_k \in \mathbb{R}^m$ to denote the center of cluster $k$. The binary variable $y_{gk}$ indicates whether group $g$ is $\alpha$-represented in cluster $k$. Let ${\cal L}_k(\textbf{z})$ be the set of allowable cluster center locations for cluster $k$ as a function of the current cluster assignments. We set ${\cal L}_k(\mathbf{z}) = \mathbbm{R}^m$ for $k$-means clustering. To represent ${\cal L}_k$ in the $k$-medians case, we introduce additional binary decision variables $d_{ik}$ to denote if data point $x^i$ is selected as the center for cluster $k$. For $k$-medians the set ${\cal L}_k(\mathbf{z})$ is defined as follows:
$${\cal L}_k(\mathbf{z})~=~\left\{c_k\in {\mathcal R}^m\;:\:c_k = \sum_{x^i \in {\cal X}} d_{ik} x^i,\quad
    d_{ik} \leq z_{ik} ~~\forall x^i \in {\cal X} ,\quad
    \sum_{x^i \in {\cal X}} d_{ik} = 1,\quad
    d_{ik} \in \{0,1\} ~~\forall x^i \in {\cal X}\right\}$$
~

We can now formulate the MR-fair clustering problem as follows:
{\footnotesize
\begin{subequations}\label{form:MRfair}\begin{align}
	\textbf{min}~~&& \sum_{x^i \in {\mathcal X}} \sum_{k \in \mathcal{K}} D(x^i,c_k) z_{ik}\label{obj:mip}\\
	\textbf{s.t.}~~&& \sum_{k \in \mathcal{K}}  z_{ik} &= 1 ~~&&\forall x^i \in {\mathcal X} \label{const:cluster_assignment}\\
    	&& \sum_{k \in \mathcal{K}} y_{gk} &\geq \beta_g &&\forall g \in {\mathcal G} \label{const:num_dominate}\\ 
	&& \sum_{x^i \in X_g} z_{ik} +  M (1 - y_{gk})&\geq \alpha \sum_{x^i \in {\mathcal X}} z_{ik}    ~~ &&\forall g \in {\mathcal G}, k \in {\mathcal K} \label{const:dominance}\\
 && u\geq \sum_{x^i \in {\mathcal X}} z_{ik} &\geq l &&\forall k \in {\mathcal K} \label{const:cardinalityA} \\
	&&z_{ik} &\in \{0,1\} &&\forall x^i \in \mathcal{X}, k \in {\mathcal K} \label{const:binary} \\
	&&y_{gk} &\in \{0,1\} &&\forall g \in {\mathcal G}, k \in {\mathcal K} \label{const:binary_y} \\
	&&c_k &\in {\cal L}_k(z) &&\forall k \in \mathcal{K} \label{const:locations}
\end{align}
\end{subequations} }
The objective \eqref{obj:mip} is to minimize the total clustering cost. Constraint \eqref{const:cluster_assignment} ensures that each data point is assigned to exactly one cluster. Constraint \eqref{const:num_dominate} enforces that each group $g$ is $\alpha$-represented in at least $\beta_g$ clusters. 
When $y_{gk} = 1$, constraint \eqref{const:dominance} ensures that group $g$ is $\alpha$-represented in cluster $k$, and when $y_{gk} = 0$,  the  constraint is trivially satisfied. The big-$M$ in   constraint \eqref{const:dominance} can be set to $\alpha n$. 
In many applications of interest, it might also be worthwhile to add a constraint on the size of the clusters to ensure that each cluster has a  minimum/maximum number of data points. Constraint \eqref{const:cardinalityA} captures this notion where $l$ and $u$ represent the lower and upper bounds for the cardinality of each cluster respectively. Note that if $u$ is finite, the big-$M$ in constraint \eqref{const:dominance} can be reduced to $\alpha u$. 
To ensure that exactly $k$ non-trivial clusters are produced by the formulation, it is important to set $l \ge 1$. Without this condition, the solution could have empty clusters where every group would be trivially $\alpha$-represented according to Definition \ref{def:alpha}.

We note that the MR-fair clustering problem is NP-hard, as the $k$-means and $k$-medians clustering problems are themselves NP-hard \citep{dasgupta2008hardness}. To solve problem \eqref{form:MRfair} in practice, we next introduce a modified version of Lloyd's algorithm called MiniReL. We also argue that, unlike other notions of fairness, first solving the clustering problem without fairness constraints and then finding a fair assignment of the data points to these cluster centers can lead to arbitrarily poor solutions. 


\subsection{MiniReL Algorithm for Fair Clustering} \label{sec:minirel_algo}
Solving the the MR-fair clustering problem \eqref{form:MRfair} to optimality is computationally challenging as it is a large scale mixed-integer optimization problem with a non-convex objective function. To solve the problem in practice, we introduce a modified version of Lloyd's algorithm which we call \textit{the Minimum Representation Fair Lloyd's Algorithm} (MiniReL) that alternates between finding cluster centers and fairly assigning data points to these centers to find a locally optimal solution, see Algorithm \ref{alg:minrepfairlloyd}. For a given (fixed) set of cluster centers $c_k$ we denote the problem \eqref{obj:mip}-\eqref{const:binary_y} the \textit{fair minimum representation assignment (FMRA) problem}, which is a linear integer program. 

\begin{algorithm}[Minimum Representation Fair Lloyd's Algorithm (MiniReL)]
\label{alg:minrepfairlloyd}

\begin{flushleft}
~\\
\textbf{Input}: Data $\mathcal{X}$, number of clusters $K$, fairness parameters $\boldsymbol\beta$, $\alpha$

\textbf{Output}: Cluster assignments $\{C_k\}_{k \in \mathcal{K}}$ and cluster centers $\{c_k\}_{k \in \mathcal{K}}$
\end{flushleft}

\begin{algorithmic}[1] \label{alg:MRfair}
\STATE Initialize $\{c_k\}_{k \in \mathcal{K}}$  (for example by uniformly at random sampling $K$ data points).
\REPEAT
\STATE Solve the FMRA \eqref{obj:mip}-\eqref{const:binary_y}  for fixed cluster centers $\{c_k\}_{k \in \mathcal{K}}$ to get optimal $\mathbf{z}$\label{alg:assign} \\
\STATE Compute the cost of solution\oo{:} \textit{current\_cost}~$=\sum_{x^i \in {\mathcal X}} \sum_{k \in \mathcal{K}} D(x^i,c_k) z_{ik}$
\FOR{$k=1,2,\dots, K$}
\IF{$k$-means}
\STATE Set  $N_k = \sum_{x^i \in \mathcal{X}} z_{ik}$, $c_k = \frac{1}{N_k} \sum_{x^i \in \mathcal{X}} z_{ik} x^i$ \label{alg:center_mean} \\
\ELSIF{$k$-medians}
\STATE Set $c_k = \argmin_{ c \in {\cal  C_k}} \sum_{x^i \in C_k} D(x^i,c)$ \label{alg:center_med}
\ENDIF
\ENDFOR
\STATE Compute the cost of new solution: \textit{improved\_cost}~$=\sum_{x^i \in {\mathcal X}} \sum_{k \in \mathcal{K}} D(x^i,c_k) z_{ik}$\label{alg:improved_cost}
\UNTIL{\textit{improved\_cost=current\_cost}}\label{alg:until} 
\STATE Set $C_k = \{x^i\in \mathcal{X} : z_{ik} = 1\} \quad \forall k \in \mathcal{K}$
\STATE \textbf{return} $\{C_k\}_{k=1}^K$ and $\{c_k\}_{k \in \mathcal{K}}$
\end{algorithmic}
\end{algorithm}

Note that the only difference between $k$-means and $k$-medians in Algorithm \ref{alg:MRfair} is how cluster centers are computed. Given a fixed set of cluster assignments (i.e., when variables $z$ are fixed in \eqref{obj:mip}-\eqref{const:binary}) the optimal choice of $c_k$ is the mean value of the data points in cluster $C_k$ for $k$-means, whereas for $k$-medians it is selected by searching across all data points in the cluster \citep{park2009simple}. In case there are multiple centers with the same cost in the $k$-medians setting we use a deterministic tie-breaking rule. For simplicity, we will refer to both variants of the algorithm as MiniReL throughout the paper. 
While the optimal assignment step in Lloyd's algorithm can be done in polynomial time, the following result shows that the FMRA problem is NP-Hard.

\begin{theorem} \label{thm:fair_assign_np}
The fair minimum representation assignment problem is NP-Hard.
\end{theorem}
\proof{Proof Sketch.} We prove that the FMRA is NP-hard via a reduction from the Exact Cover by 3-Sets (X3C) problem. Given an instance of X3C, we first construct a bipartite graph $G = (U\cup W, E)$, where nodes in $U$ correspond to elements, and nodes in $W$ to candidate sets. 
Using this graph, we construct an FMRA instance which has a data point for each node and points associated with nodes in $U$ and $W$ belong to different groups. We then use MR-fairness constraints to ensure that the clustering corresponds to a solution to the X3C problem, if one exists. For a full proof see Appendix \ref{app:fair_assign_np}. 

Note that if the FMRA problem is infeasible for any given set of cluster centers, it cerifies that no MR-fair clustering exists with the given $\alpha$ and $\boldsymbol\beta$. While integer programs do not always scale well to large datasets, our computational experiments show that FMRA can be solved to optimality in a reasonable amount of time even for datasets with tens of thousands of data points. In Section \ref{sec:scaling} we describe the computational techniques that help scale our algorithm. 

A natural question is whether MiniReL converges to a locally optimal solution as Lloyd's algorithm does. When discussing local optimality, it is important to formally define the local neighborhood of a solution. In the absence of fairness constraints, data points must be assigned to the closest centers to minimize cost. Consequently, a clustering is locally optimal if perturbing the centers does not improve the clustering cost.
In our setting, we define a local change as any perturbation to a cluster center, \textit{or} an individual change to cluster assignment (i.e., moving a data point from one cluster to another). With this notion of a local neighborhood, we next show that the MiniReL algorithm converges to a local optimum in finite time. Note that while MiniReL converges to a local optimum, the solution may be arbitrarily worse than the global optimum as is the case for Lloyd's algorithm.

\begin{theorem} \label{thm:finite_convergence} 
MiniReL converges to a local optimum in finite time.
\end{theorem}
\proof{Proof.} See Appendix \ref{app:pf_convergence}.

\subsection{Checking Feasibility of Fairness Parameters} \label{sec:feasibility}
\newcommand{\combo}{t}
\newcommand{\allcombo}{\mathcal T} 


The MR-fair clustering problem is not guaranteed to be feasible for every combination of fairness parameters ($\alpha, \boldsymbol{\beta}$) and cardinality constraints ($ l, u$). 
In this section, we construct a {small} integer program that can check whether a set of fairness parameters and cardinality constraints are feasible for a given dataset. We start with enumerating every possible combination of sensitive features present in the dataset, which we call a \textit{type}. For instance, consider a dataset with two sensitive features: Gender (Male, Female, Nonbinary) and Age (Youth, Adult). The possible types are (Male, Youth), (Female, Youth), (Non-binary, Youth), (Male, Adult), (Female, Adult), (Non-binary, Adult). 
Let  ${\allcombo}$ represent the set of all possible types, ${\allcombo}_g$ be the set of combinations that involve a group $g \in {\cal G}$, and $n_\combo$ represent the number of data points in ${\cal X}$ that have the type $\combo\in{\allcombo}$.
We next formulate the \textit{MR-fair feasibility problem} as an integer program where variable $z_{\combo k}$ denotes the number of data points of type $\combo\in\allcombo$ assigned to cluster $k$ and,  variables $s_l,s_u,s_g$ keep track of constraint violations:

 
{\footnotesize
\begin{subequations}\label{from:MRfeasible}\begin{align}
	\textbf{min}~~&& \lambda_l s_l + \lambda_u s_u&+\sum_{g \in {\mathcal G}} \lambda_gs_g\label{obj:feasibility_check}\\
	\textbf{s.t.}~~&& \sum_{k \in \mathcal{K}}  z_{\combo k} &= n_\combo  ~~&&\forall \combo  \in {\allcombo} \label{const:feas_cluster_assignment}\\
	&& \sum_{\combo  \in {\allcombo}_g} z_{\combo k} +  M (1 - y_{gk})&\geq \alpha \sum_{\combo  \in {\allcombo}} z_{\combo k}    ~~ &&\forall g \in {\mathcal G}, k \in {\mathcal K} \label{const:feas_dominance}\\
	&& \sum_{k \in \mathcal{K}} y_{gk} &\geq \beta_g - s_g &&\forall g \in {\mathcal G} \label{const:feas_num_dominate}\\ 
 && u + s_u\geq \sum_{\combo  \in {\allcombo}} z_{\combo k} &\geq l - s_l &&\forall k \in {\mathcal K} \label{const:cardinalityC} \\
	&&z_{\combo k} &\in \mathbb{Z}_+ &&\forall \combo  \in \allcombo, k \in {\mathcal K} \label{const:feas_binary} \\
	&&y_{gk} &\in \{0,1\} &&\forall g \in {\mathcal G}, k \in {\mathcal K} \label{const:feas_binary_y} \\
    && s_g &\in \mathbb{Z}_+ &&\forall g \in {\mathcal G}
\end{align}\end{subequations}
}
The objective \eqref{obj:feasibility_check} minimizes the (weighted) parameter adjustment needed to make the MR-fair clustering problem feasible. If the optimal objective value is $0$, then the original problem is feasible - otherwise the solution identifies the necessary changes to the fairness and cardinality parameters to make the problem feasible. Constraints \eqref{const:feas_cluster_assignment}-\eqref{const:feas_binary_y} mirror the constraints from the FMRA problem, except the constraints now operate on combinations $\combo$ instead of individual data points. 
If the problem is feasible, the resulting clustering can be disaggregated to construct a feasible solution to the MR-fair clustering problem.

Note that the MR-fair feasibility problem \eqref{from:MRfeasible} has only $|{\cal G}|K$ binary variables and $|{\allcombo}|K$ integer variables. 
In a practical setting, $|{\cal G}|$ and $|{\allcombo}|$ are both small numbers, and are typically orders of magnitude smaller than the number of data points $n$. For example in our experiments $n$ can be up to $50,000$ while $|{\cal G}| \leq 5$ and $|{\allcombo}| \leq 6$. This means that the MR-fair feasibility problem can be solved quickly in practice either to verify the feasibility of the input parameters, or to obtain minimal changes to them to ensure feasibility. Moreover, when $|{\allcombo}|$, $|{\cal G}|$, and $K$ are treated as fixed parameters, the MR-fair feasibility problem can be solved in polynomial time.

\begin{proposition}
    For fixed $|{ \allcombo}|$, and $K$, the MR-fair feasibility problem \eqref{from:MRfeasible} can be solved in polynomial-time.
\end{proposition}
\proof{Proof.}
    Follows from the fact that integer programs can be solved in polynomial-time when the number of variables (dimension) is fixed. See \cite{lenstra1983integer,eisenbrand2003fast}.
\endproof

\subsection{A Natural Solution Approach and an Inapproximability Result} \label{sec:inapprox}
One natural approach \citep{bera2019,esmaeili2021fair,esmaeili2022fair} used for other notions of fairness in clustering is to first obtain cluster centers by solving the clustering problem without fairness constraints and then find a fair assignment of the data points to these centers. In some settings, this has been shown to provide solutions with a provable approximation guarantee even when the centers are only approximately optimal for the clustering problem without fairness constraints \citep{bera2019,esmaeili2021fair,esmaeili2022fair}. 
Unfortunately, we show that in our MR-fairness setting this approach can lead to arbitrarily bad solutions. Let $z^*, c^*$ be the optimal solution to the MR-fair clustering problem \eqref{form:MRfair}. Let $c^*_{UF}$ be the optimal centers for the clustering problem without fairness constraints and $z^*_{FA}$ be the optimal (fair) assignment of data points to centers $c^*_{UF}$. Note both $z^*$ and $z^*_{FA}$ are feasible fair assignments. Let $COST(z, c) = \sum_{x^i \in {\mathcal X}} \sum_{k \in \mathcal{K}} D(x^i, c_k) z_{ik}$ 
denote the cost of the clustering.

\begin{theorem} \label{thm:fair_adjustment}
There does not exist a constant $\Gamma > 0$ such that:
$$
COST(z^*_{FA}, c^*_{UF}) \leq \Gamma \cdot  COST(z^*, c^*)
$$
In other words, fairly assigning data points to (approximately) optimal unfair centers can lead to arbitrarily worse performance relative to the optimal solution of the MR-fair clustering problem.
\end{theorem}

\begin{figure}[t]
    \centering
    \begin{subfigure}
      \centering
      \includegraphics[width=0.65\textwidth]{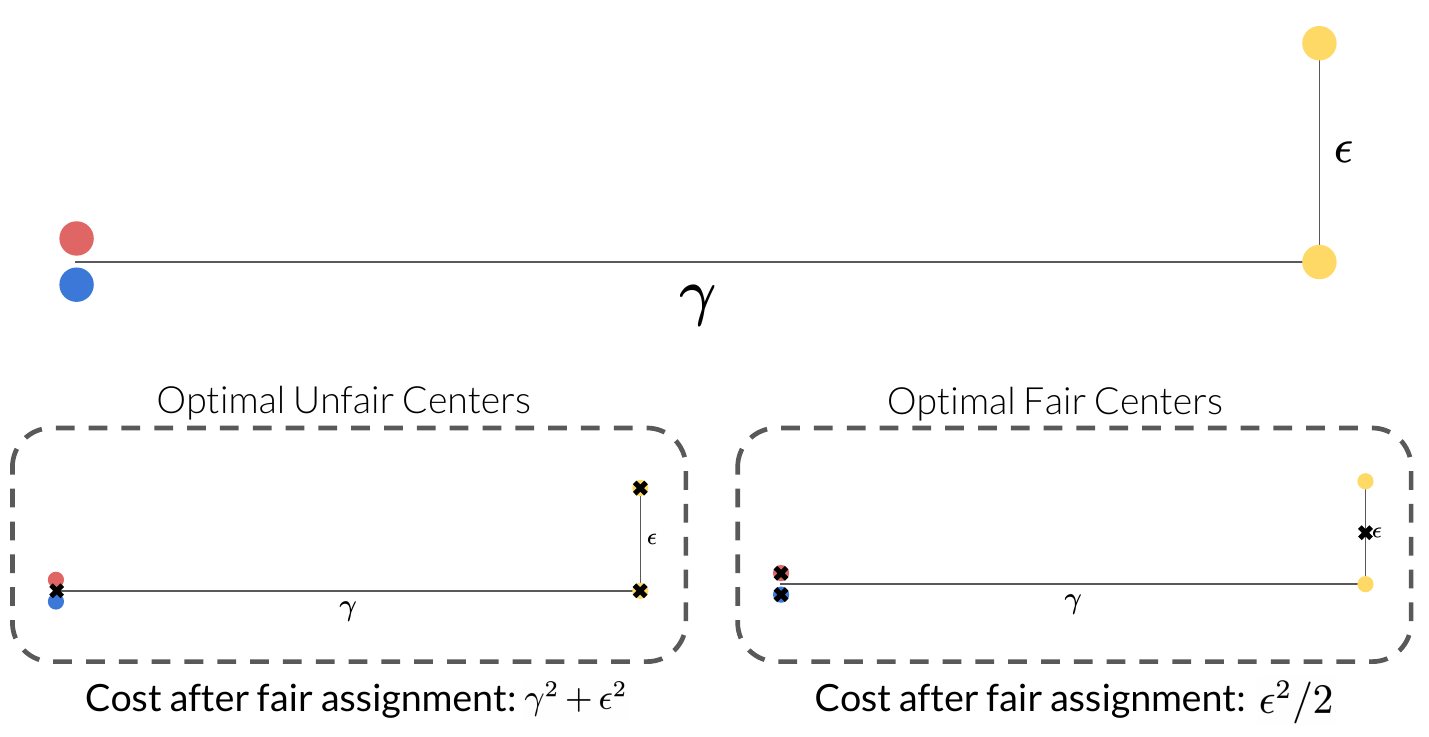}
    \end{subfigure}\vskip-.1cm 
  \caption {\label{fig:inapprox} (Top) Bad instance of MR-fair clustering problem with four data points and three groups (represnted by colors). (Bottom Left) Optimal centers for the unfair clustering prioblem. (Bottom Right) Optimal centers for the MR-fair clustering problem.
  \vspace{-0.5cm}
  }
\end{figure}

\proof{Proof.}
To prove the claim we will construct an instance of the MR-fair $k$-means problem in two dimensions with four data points and three groups (Red, Blue, and Yellow). The first two points are located at $(0,0)$ and belong to the red and blue groups, respectively. The next two data points belong to the yellow group and are located at $(\gamma,0)$ and $(\gamma,\epsilon)$ for some $\gamma,\epsilon>0$ (see Figure \ref{fig:inapprox}).
We set the parameters for the MR-fair problem as follows: $K=3$, $\alpha > 0.5$, and $\beta_{R}=\beta_{Y}=\beta_{B}=1$. 
For the remainder of this proof we use $D(x^i, c_k)= \lVert x^i- c_k\lVert^2_2$  but the proof goes through with any other distance metric.


Note that any optimal solution to the clustering problem without fairness constraints will place the cluster centers $c^*_{UF}$  at $(0,0),$ and $(\gamma,0), (\gamma,\epsilon)$ yielding a total clustering cost of $0$. Now consider a fair assignment of points to these fixed centers. It is straightforward to see that, under the fairness constraints, one must assign either the red or blue data point to the center at $(\gamma, 0)$ and assign the yellow data point at $(\gamma, 0)$ to the center at $(\gamma, \epsilon)$ yielding a cost of $\gamma^2 + \epsilon^2$. 

Now consider the optimal solution to the fair clustering problem which places the cluster centers at $(0,0), (0,0), (\gamma,\epsilon/2)$. The red and the blue data point are assigned to one of the centers at $(0,0)$, and both yellow data points are assigned to the center at $(\gamma,\epsilon/2)$ which satisfies the fairness constraints and incurs a total cost of $\epsilon^2/2$. The ratio of the costs of the two solutions is $\frac{\gamma^2 + \epsilon^2}{\epsilon^2/2}$ which can be made arbitrarily large by increasing $\gamma$, completing the proof. Note that for $k$-medians the same construction, with optimal fair centers at $(0,0), (0,0), (\gamma,0),$ proves the same result. We also note that this proof holds even when there is a small distance (less than $\epsilon$) between the red and blue points.\qed
\endproof

Theorem \ref{thm:fair_adjustment} shows that the cost of the optimal fair clustering can be arbitrarily larger than that of the optimal clustering without fairness. We define \textit{the price of fairness} (PoF) to be the ratio of the fair clustering cost to the clustering cost without fairness:
$$
PoF = \frac{\text{Clustering Cost with Fairness Constraint}}{\text{Clustering Cost of Agnostic Solution}}
$$

Using this definition Theorem \ref{thm:fair_adjustment} implies the following:

\begin{corollary}
For any $k \geq 2$, imposing MR-fairness can lead to an unbounded Price of Fairness.
\end{corollary}


		

\section{Scaling MiniReL: Two-Stage Decomposition} \label{sec:scaling}

The main computational bottleneck in the MiniReL algorithm is solving the FMRA problem, which jointly determines which groups are $\alpha$-represented in which clusters (i.e., setting $y_{gk}$ variables) and assigns data points to clusters to satisfy MR-fairness constraints  (i.e., setting $z_{ik}$ variables). 
This is a computationally demanding problem due to its scale (there is a $K$ binary variables for each data point), as well as the inherent symmetry in the problem formulation, and the use of big-M constraints (that lead to weak linear relaxations \citep{conforti2014integer}). To address these challenges, we introduce a heuristic two-stage decomposition scheme that decouples the selection of the $y$ and $z$ variables into two sequential problems. 

The first stage problem, which we call the \textit{\firststage (\fsabbrvx)}, determines which groups are $\alpha$-represented in which clusters (the $y_{gk}$ variables). 
The second stage problem, which we call the \textit{\secondstage (\ssabbrvx)}, assigns data points to clusters consistent with the first stage decisions (the $z_{ik}$ variables).  One approach to solve the first-stage problem \fsabbrv is to relax the FMRA problem by allowing $z$ variables to take fractional values (i.e., problem \eqref{obj:mip}-\eqref{const:locations} with fixed centers, $y_{gk}$ binary and $z_{ik} \in [0,1]$). This significantly reduces the number of binary variables, making the problem much faster to solve than the FMRA problem in practice. While this approach does not guarantee a feasible second-stage solution, in Section \ref{sec:flow}, we show that the fractional first stage solution can be rounded to an integer solution with only small violation to the fairness constraints. In Section \ref{sec:prefix} we also present a heuristic to solve the first-stage problem. 

As the $y_{gk}$ variables are fixed in the second stage problem (\ssabbrvx), the  big-M constraints \eqref{const:feas_dominance} are not needed. Furthermore, this also breaks symmetry in the IP (i.e., removes permutations of feasible cluster assignments), dramatically improving the computation time. However, despite this simplification, the following result shows that the second-stage problem is still NP-Hard.

\begin{theorem} \label{thm:second_stage_np}
The \secondstage (\ssabbrvx) is NP-Hard.
\end{theorem}

\proof{Proof Sketch} We prove that the \secondstage is NP-hard via a reduction from the 3-SAT problem. The key idea behind the reduction is to construct a corresponding FMRA instance with two clusters and one data point for each variable in the 3-SAT instance. One cluster corresponds to setting that variable to be True, and the other to False. MR-fairness constraints ensure that any solution to the FMRA instance corresponds to a truth assignment that satisfies all clauses in the original 3-SAT instance. For a full proof see Appendix \ref{app:second_stage_np}. \qed

Despite this negative result, in Section \ref{sec:flow} we present a polynomial-time algorithm that solves \ssabbrv  with a small additive fairness violation. 
We emphasize that this two-stage decomposition scheme is not guaranteed to solve the FMRA problem to optimality. However, it performs well in practice, finding near-optimal solutions very quickly (see Section \ref{sec:algo_design}).
While  both stages can be solved at every iteration of the MiniReL algorithm, we observed that the first stage solution often remained unchanged. To further improve the computation time, we also present a \textit{pre-fixed} version of the algorithm where we only set the $y_{gk}$ variables once at the beginning of the algorithm. In problems where only a single group can have $\alpha$-representation in a cluster (i.e., a data point can only be part of one group and $\alpha > 0.5$), pre-fixing preserves an optimal solution to the FMRA problem. However, in more complicated settings this is not the case. It is worth noting that the MiniReL algorithm is itself a heuristic, and thus the pre-fixing scheme has ambiguous effects on the cost of the solution as it may cause the algorithm to converge to a better local optimum. Figure \ref{fig:decomp} illustrates the three approaches to solve the overall problem.

In the rest of this section, we describe polynomial-time heuristic algorithms to solve the first stage ~(\fsabbrvx) and the second stage (\ssabbrvx) problems.

\begin{figure}[t]
    \centering
    \begin{subfigure}
      \centering
      \includegraphics[width=0.8\textwidth]{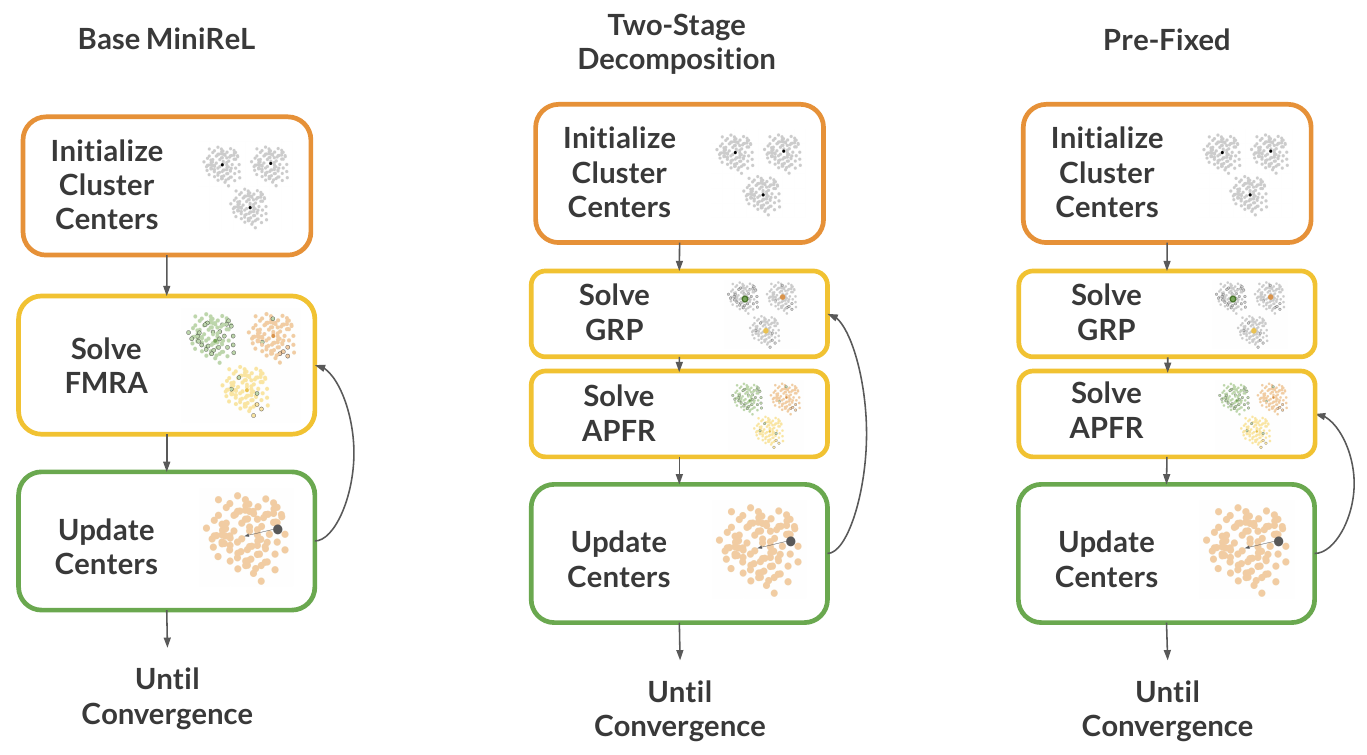}
    \end{subfigure}\vskip-.31cm 
  \caption {\label{fig:decomp} (Left) Base MiniReL algorithm where the assignment is doing by solving the full FMRA problem. (Center) MiniReL Algorithm with two-stage decomposition approach for solving the FMRA problem where we sequentially set $y_{gk}$ variables and $z_{ik}$ variables. (Right) Pre-fixed approach where the first stage problem to set $y_{gk}$ variables is only done once at beginning of algorithm.
  \vspace{-0.5cm}
  }
\end{figure}

\subsection{Solving \fsabbrv via Aggregation}\label{sec:prefix}

We next present a polynomial-time heuristic to solve the \fsabbrv problem. The starting point of our algorithm is a given clustering of the data points together with the corresponding cluster centers. This initial clustering is not required to be MR-fair, and can generated by first running Lloyd's algorithm withour fairness constraints. In Appendix \ref{app:warm_start} we show experiments with different initial starting clusters.

For a cluster $k$ and group $g$, let $q_{kg}\ge0$ be the additional number of points from group $g$ needed to make this group $\alpha$-represented in cluster $k$ (i.e., smallest integer $q_{kg}\ge0$ that satisfies $q_{kg} + |C_k\cap X_g| \geq \alpha(q_{kg} + |C_k|)$). Let $c_{(x)} = \argmin_{ c \in \{c_1, \dots, c_K\}} \{\|x - c\|_2^2\}$ denote the closest center for point $x \in {\cal X}$. 
We estimate the (myopic) cost increase to make group $g$ $\alpha$-represented in cluster $k$ as follows:
$$
    m_{gk} = \min_{S \subset X_g \setminus C_k: |X| = q_{kg}} \sum_{x \in S} \big(D(x, c_k) - D(x, c_{(x)}))\big)
    $$
 Note that this does not fully capture the effect of moving the data points in $S$ to cluster $k$ within the overall problem. We can now formulate the problem of performing the group representation assignment as follows:

{\footnotesize
\begin{subequations}\label{form:prefix}\begin{align}
	\textbf{min}~~&& \sum_{(g,k) \in {\mathcal G} \times {\mathcal K}} m_{gk} y_{gk}\label{obj:prefix}\\
	\textbf{s.t.}~~&& \sum_{k \in \mathcal{K}}  y_{gk} &\geq \beta_g ~~&&\forall g \in {\mathcal G} \label{const:beta_const}\\
	&& \sum_{k \in \mathcal{K}}  z_{\combo k} &= n_\combo ~~&&\forall \combo \in {\mathcal \allcombo } \label{const:prefix_cluster_assignment}\\
	&& \sum_{\combo \in {\cal \allcombo }_g} z_{\combo k} +  M (1 - y_{gk})&\geq \alpha \sum_{\combo \in {\cal \allcombo }} z_{\combo k}    ~~ &&\forall g \in {\mathcal G}, k \in {\mathcal K} \label{const:prefix_dominance}\\
 && u \geq \sum_{\combo \in {\mathcal \allcombo }} z_{\combo k} &\geq l &&\forall k \in {\mathcal K} \label{const:prefix_cardinality} \\
	&&z_{\combo k} &\in \mathbb{Z}_+ &&\forall \combo  \in \mathcal{\allcombo }, k \in {\mathcal K} \label{const:prefix_binary} \\
	&& y_{gk} &\in \{0,1\} ~~ &&\forall g \in \mathcal{G},k \in \mathcal{K} \label{const:prefix_binary}
\end{align}\end{subequations}
}
The variables and constraints in this formulation are similar to that of the MR-fair feasibility problem \eqref{from:MRfeasible}, except we do not have the $s$ variables.
The objective \eqref{obj:prefix} is a proxy for the cost of pre-fixing. 
Similar to the MR-fair feasibility problem, for fixed $|{\cal K}|$, $|{\cal T}|$ the prefix problem can be solved in polynomial time. In practice, this is a very small IP that can be solved in seconds. 

The key computational benefit of the heuristic is is due to the fact that it aggregates the data points with the same combination of sensitive features. However, a drawback of this aggregation is that it may fail to correctly estimate the true clustering cost. For instance, the objective term may `double-count' individual data points and underestimate the true objective. In initial experiments we observed that the objective value of this heuristic could underestimate the objective of the first iteration of the APFR by up to 20\% which lead to worse clustering performance. To mitigate this issue we use a \textit{fallback mechanism} that solves the GRP via an integer program if the objective of the prefix heuristic differs from the objective of the APFR by more than 10\%. While this feedback was employed in less than 10\% of the instances evaluated in our experiments, it reduced the worst-case increase in clustering cost (relative to the baseline MiniReL algorithm) from over 20\% to under 2\%.

\subsection{Solving \ssabbrv via Network Flow Rounding} \label{sec:flow}
In this section, we present a polynomial-time algorithm to solve the second-stage \ssabbrv problem while approximately satisfying the fairness constraints. The main idea is to take the fractional assignment of data points obtained coming from the first stage \fsabbrv problem (or from solving the linear relaxation of APFR) and round it to an integer assignment using a min-cost network flow model. Recall that in the \ssabbrv problem, clusters have already been assigned groups that must be $\alpha$-represented in them (i.e., the $y_{gk}$ variables are fixed). We start by introducing a graph formulation for the network flow rounding problem, before explaining how we use it to solve the \ssabbrv problem and analyzing its properties.

\paragraph{Network Flow Graph:} Given a fractional solution $\mb{z}^{LP}$, we construct a flow graph $G=(V,E)$ where $V$ consists of three sets of nodes $V_x$, $V_{TK}$, and $V_K$. For every data point $x^i \in {\cal X}$ we construct a node $v_i \in V_{X}$ with a supply of 1. For every cluster $k$ we create a set of nodes $V_{TK}$ where we have 
one node for every possible 
\textit{type} associated with groups that must be $\alpha$-represented in the cluster 
plus one `remainder' node  for all other types (see Section \ref{sec:feasibility} for an expanded discussion and example of types). 
Let $T_k$ be the set of types for a cluster $k$ including the remainder node for all other types. For each $t \in T_k$ we create create one node $v_{tk} \in V_{TK}$ with demand $d_{tk} = \lfloor \sum_{x^i \in \bigcap_{g \in t} X_g} z_{ik}^{LP} \rfloor$. We create an edge with capacity one between data point $x^i$ and every node $v_{tk}$ if $x^i \in \cap_{g \in t} X_g$. The cost for these edges is equal to the clustering cost of assigning $x^i$ to the fixed center $c_k$ (i.e., $D(x^i, c_k)$). The costs of all other edges in this graph are 0.

\begin{figure}[t]
    \centering
    \begin{subfigure}
      \centering
      \includegraphics[width=0.5\textwidth]{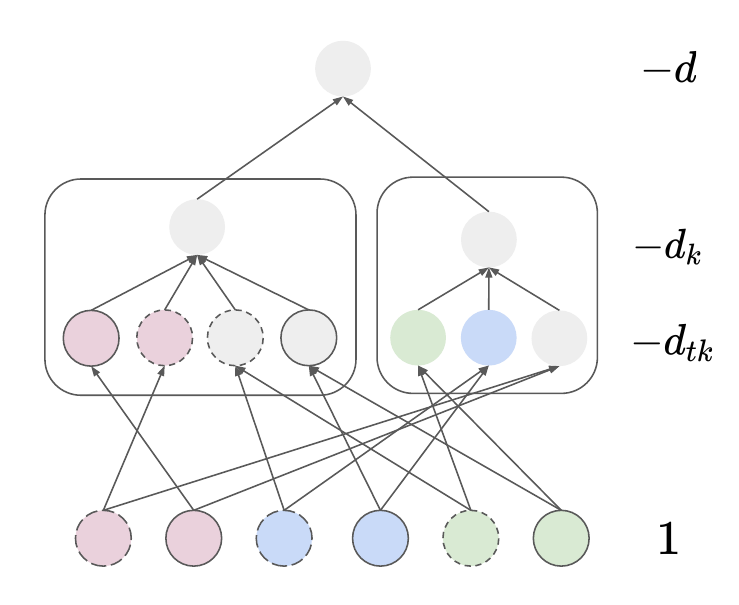}
    \end{subfigure}\vskip-.31cm 
  \caption {\label{fig:flow_fmra} Sample network flow formulation for rounding the FMRA problem under pre-fixing with six data points, two clusters, and two sensitive attributes: Color (Purple, Blue or Green) and Border Type (Dashed or Solid). After pre-fixing, cluster 1 has to $\alpha$-represent both the Red and Dashed groups, and cluster 2 has to represent the green and blue groups. Node supplies (positive) and demands (negative) are indicated on the right hand side. 
  \vspace{-0.5cm}
  }
\end{figure}

For every cluster $k$ we also create a node $v_k \in V_{K}$ with demand $d_k = \lfloor \sum_{x^i \in {\cal X}} z_{ik}^{LP} \rfloor - \sum_{t \in T_k} d_{tk}$. We create an edge with capacity one between every node $v_{k}$ and every node $v_{tk}$ that belong to the same cluster. Finally, we add a sink note $v_{s}$ with demand $d= |{\cal X}| - \sum_{k} (d_k + \sum_{t \in T_k} d_{tk})$ to capture any residual supply. For every node $v_k \in V_k$ with fractional weight assigned to it in the LP solution (i.e., $\sum_{x^i \in {\cal X}}z^{LP}_{ik} \notin \mathbb{Z}$) we create an edge with capacity one between it and the node $v_t$.

\paragraph{Example} Consider a simple example of the \ssabbrvx~with 2 clusters, 6 data points, and two sensitive features: color (red, blue, and green) and border type (dotted and solid). The six data points correspond to every unique combination of sensitive features. The first cluster needs to be $\alpha$-represented by the red and dotted-border group, and teh second by the blue and green groups. We now construct the flow graph corresponding to this instance (see Figure \ref{fig:flow_fmra}). Each data point is given an associated node and supply of one. The first cluster has nodes $V_{t1}$ corresponding to types: (red, solid), (red, dotted), (not red, dotted), (not red, not dotted). The second cluster has nodes $V_{t2}$ correspond to types: (green, any), (blue, any), (not green or blue, any). Each of the nodes in $V_{t1}, V_{t2}$ has demand equal to $d_{tk}$. Each cluster also includes a node $v_k$ with demand equal to $d_k$. Finally we include a sink node $v_s$ with demand $d$. 

\paragraph{Flow-APFR} To construct a solution to the \ssabbrvx~problem using the network flow model we solve a minimum-cost flow problem over the graph. Note that this min-cost flow problem can be solved in polynomial time via standard min-cost flow algorithms \citep{williamson2019network}. We denote the minimum cost flow problem associated with this graph the \textit{Flow-\ssabbrvx}~problem. We now show that by solving the Flow-\ssabbrvx~problem we can generate a good assignment of data points to clusters:

\begin{theorem} \label{thm:flowfmra_binary}
Solving the Flow-\ssabbrvx~problem yields a binary assignment of data points to clusters that satisfies the cardinality constraints and has a cost that does not exceed the cost of the fractional solution $\mb{z}^{LP}$.
\end{theorem}

\proof{Proof} Recall that $\mb{z}^{LP}$ is the optimal fractional assignment of the \ssabbrvx~problem (found by solving the LP relaxation or as the output of the \fsabbrvx~problem). Note that since all capacities and demands in the Flow-\ssabbrvx~problem are integer, the optimal solution is an integer flow. We interpret this flow as a binary assignment of data points to clusters by setting $\bar{z}_{ik} = 1$ if there is a unit of flow between node $v_i \in V_{\cal X}$ and any node $v \in V_{Tk}$, and 0 otherwise. By construction every data point will be assigned to exactly one cluster. Since the LP solution is feasible for the network flow problem, we know that the cost of the solution $\bar{\textbf{z}}$ is lower or equal to the cost of the LP solution.

Note that at least $\lfloor \sum_{x^i \in {\cal X}} z_{ik}^{LP} \rfloor$ and at most $\lceil \sum_{x^i \in {\cal X}} z_{ik}^{LP} \rceil$ data points are assigned to cluster $k$. This follows from the fact that a total of $d_k + \sum_{t \in T_k} d_{tk} = \lfloor \sum_{x^i \in {\cal X}} z_{ik}^{LP} \rfloor$ demand is consumed in nodes corresponding to cluster $k$ and the edge capacity of the edge of 1 between $v_k$ and $v_t$ ensures at most one additional point is allotted to cluster $k$. We now claim that the solution to the Flow-\ssabbrvx~problem satisfies the cardinality constraints. Suppose it did not, then $\lfloor \sum_{x^i \in {\cal X}} z_{ik}^{LP} \rfloor < l$ or $\lceil \sum_{x^i \in {\cal X}} z_{ik}^{LP} \rceil > u$. This implies that $\sum_{x^i \in {\cal X}} z_{ik}^{LP} < l$ or $\sum_{x^i \in {\cal X}} z_{ik}^{LP} > u$, which would violate constraints \eqref{const:cardinalityA} and contradicts $\mb{z}^{LP}$ being a solution to the LP relaxation of the \ssabbrvx~problem. \qed

It remains to analyze the impact of this rounding scheme on the fairness of the final solution. In particular, we look at the impact of this rounding procedure on the $\alpha$-representation fairness constraints (i.e., constraints \eqref{const:dominance} for each group $g$ and cluster $k$ where $y_{gk}=1$). Recall the following definition for an additive constraint violation:


\begin{definition}[Additive Constraint Violation]
For a given constraint $g(x) \leq 0$, a solution $\bar{x}$ is said to have an additive constraint violation of $\delta$ for $\delta \geq 0$ if
$g(\bar{x}) \leq  \delta$.
\end{definition}

We now show that the assignment produced by the Flow-\ssabbrvx~problem has a small additive violation of the MR-fairness constraints. 

\begin{theorem} \label{thm:flowfmra}
Let $\gamma = \min(\lceil \alpha^{-1} \rceil, \max_{f \in {\cal F}}|{\cal G}_f|)$. The additive fairness violation of the clustering produced by the Flow-\ssabbrvx~problem is at most:

$$
\gamma^{|{\cal F}| - 1} + \alpha \mathbb{I}(\gamma > 2)
$$
\end{theorem}

\proof{Proof} 
Let $\bar{\mb{z}}$ be the binary assignment  produced by the Flow-\ssabbrvx~problem and consider the additive fairness violation of the constraint for group $g$ in cluster $k$ by $\bar{\mb{z}}$. We can write this violation as:
$$
\delta = \alpha\sum_{x^i \in {\cal X}}\bar{z}_{ik} - \sum_{x^i \in X_g}\bar{z}_{ik}  =\alpha \sum_{x^i \notin X_g} \bar{z}_{ik}  - (1-\alpha)\sum_{x^i \in X_g} \bar{z}_{ik}.
$$
Let $f_g, f_{g'} \geq 0$ be the difference between the fractional and integer assignment for the group $g$ and data points outside $g$ respectively:
$$
    f_g = \sum_{x^i \in X_g} z_{ik}^{LP} - \sum_{x^i \in X_g} \bar{z}_{ik} \quad\quad
    f_{g'} = \sum_{x^i \notin X_g} \bar{z}_{ik}  - \sum_{x^i \notin X_g} z_{ik}^{LP}
$$

Rewriting the fairness constraint violation we get:{\small\begin{align*}
\delta 
&= \alpha \sum_{x^i \notin X_g} z_{ik}^{LP} + \alpha f_{g'} -
(1-\alpha)\sum_{x^i \in X_g} z_{ik}^{LP} - (1-\alpha)f_g  \\
&= (1-\alpha)f_{g} + \alpha f_{g'} + ( \alpha \sum_{x^i \notin X_g} z_{ik}^{LP} -  (1- \alpha)\sum_{x^i \in X_g} z_{ik}^{LP}) \\
&\leq (1-\alpha)f_{g} + \alpha f_{g'}
\end{align*}
}
where the final inequality is implied by the feasibility of the LP solution for the fairness constraint. Let ${\cal F}_k$ be the set of sensitive features $f$ with at least one group $g \in {\cal G}_f$ $\alpha$-represented in cluster $k$, and let $\eta_{fk} = \max(1+\sum_{g \in {\cal G}_f} y_{gk}, |{\cal G}_f|)$. Note that there are a total of $\prod_{f \in {\cal F}_k} \eta_{fk}$ nodes in each cluster $k$. Let $V_{gk} \subset V_{TK}$ be the set of nodes $v_{tk}$ such that $g \in t$, and note $|V_{gk}| = \prod_{f \in {\cal F}_k: \\ g \notin {\cal G}_f} \eta_{fk}$. Also note that the $V_{gk}$ sets do not form a partition of $V_{TK}$ as a node $v_{tk}$ can correspond to multiple groups. Similarly let $V_{g'k}$ be the set of nodes $v_{tk}$ such that $g \notin t$. Note that $|V_{gk}| \leq| V_{g'k}|$, with the comparison being strict only when $|{\cal G}_f| > 2$. By construction we know that for every node $v_{tk} \in V_{gk}$ at least $d_{tk} = \lfloor \sum_{x^i \in \cap_{g \in c} X_g} z_{ik}^{LP} \rfloor $ units of flow are routed to it. Thus $f_g \leq \sum_{v_{tk} \in V_{gk}} \sum_{x^i \in \cap_{g \in c} X_g} z_{ik}^{LP} - d_{tk} \leq |V_{gk}|$. By a similar argument we also have $f_{g'} \leq |V_{g'k}|$. By construction the total number of data points assigned to cluster $k$ is at most $\lceil \sum_{x^i \in {\cal X}} z_{ik}^{LP} \rceil$, which means $f_{g'} - f_{g} \leq 1$. Therefore the worst-case fairness violation from rounding can be seen as the following optimization problem:
$$
\max_{f_g, f_{g'} \geq 0} (1-\alpha)f_{g} + \alpha f_{g'} \quad \textbf{s.t.}~~~~ f_{g'} - f_{g} \leq 1, ~~ f_{g'} \leq |V_{g'k}|,~~ f_{g} \leq |V_{gk}|
$$
Recall that $|V_{gk}| \leq |V_{g'k}|$ and thus by inspection, the optimal solution is $(1-\alpha)|V_{gk}| + \alpha \min(1+|V_{gk}|, |V_{g'k}|)$. Note that $|V_{g'k}| > |V_{gk}|$ only if $\max_{f \in {\cal F}} |{\cal G}_f| > 2$ AND $\lceil \alpha^{-1} \rceil > 2$ (i.e., there can be $\geq 3$ nodes corresponding to the same sensitive feature). This condition is equivalent to saying $\gamma > 2$. Finally, we complete the proof by bounding $|V_{gk}|$. Recall that by construction: $$
|V_{gk}| ~= \prod_{f \in {\cal F}_k: g \notin {\cal G}_f} \eta_{fk} ~\leq \prod_{f \in {\cal F}_k: g \notin {\cal G}_f} \gamma ~\leq ~\gamma^{|{\cal F}|-1}
$$\qed


We now show that this result generalizes the earlier result of \citep{bercea2018cost} and gives an additive fairness violation of 1 in the special case of two disjoint groups.

\begin{corollary} \label{corr:flowfmra_twogroup}
    In the special case when $|{\cal G}| = 2$ and $|{\cal F}|= 1$, the FMRA-flow problem guarantees of an additive fairness violation of at most 1.
\end{corollary}

\proof{Proof of Corollary \ref{corr:flowfmra_twogroup}}
Follows from the fact that in the two group case $|{\cal F}| = 1$ and $\gamma = 2$.

One limitation of Theorem \ref{thm:flowfmra} is the the worst-case fairness violation scales exponentially in the number of sensitive features $|{\cal F}|$. In the following result, we show that the fairness violation is also upper bounded by $K$ and the number of $\alpha$-representation constraints and thus scales linearly with $|{\cal F}|$:

\begin{theorem}
The additive fairness violation of the clustering produced by the Flow-\ssabbrvx~problem is at most:
$$
K + \sum_{g \in {\cal G}}\beta_g  
\leq K + K \gamma |{\cal F}|
$$	
\end{theorem}
\proof{Proof} We prove this result via a counting argument based on the linear relaxation of the \ssabbrvx~problem. 
Let $\mb{z}^{LP}$ be the optimal solution to this linear relaxation. Take all $z_{ik}^{LP}$ variables with values in $\{0,1\}$ and fix them, then re-solve the LP. We now have a LP where all basic $z_{ik}^{LP}$ variables are fractional. We now bound the number of data points with fractional variables, denoted $\hat{{\cal X}}$. Let $n_{frac}$ be the number of fractional basic $z_{ik}^{LP}$ variables. Note that for every data point $x^i \in \hat{\cal X}$ there must be at least two basic fractional variables $z_{ik}^{LP}$ associated with it, implying $2 |\hat{{\cal X}}| \leq n_{frac}$. 

Also note that the LP has at most the following tight constraints:
\begin{itemize}
\item $|\hat{{\cal X}}|$ constraints corresponding to constraint \eqref{const:cluster_assignment}.
\item $\sum_{g \in {\cal G}} \beta_g \leq K \gamma |{\cal F}|$ constraints of type \eqref{const:dominance} corresponding to the pre-fixed $\alpha$-representation constraints.
\item $K$ active upper or lower bound constraints of type \eqref{const:cardinalityA}.
\end{itemize}
In total there are at most $K + K \gamma_{MAX} + \hat{{\cal X}}$ tight constraints, implying $n_{frac} \leq  K + K \gamma_{MAX} + \hat{{\cal X}}$.
Combining both the upper and lower bound for $n_{frac}$ and re-arranging terms we get that $\hat{{\cal X}} \leq K + \sum_{g \in {\cal G}} \beta_g$. The worst-case fairness violation is upper-bounded by the number of fractional variables, completing the result. \qed

\section{Numerical Results} \label{sec:exp}

To benchmark our approach, we evaluate different variants of MiniReL on a suite of 12 datasets from the UCI machine learning repository \citep{dua2017uci} that have been previously used in the fair clustering literature. Table \ref{tab:big_table} summarizes the size of each dataset with respect to the number of data points $n$, number of features $m$, the number of sensitive features $|{\cal F}|$, and the total number of groups $|{\cal G}|$. We also include results for the Brunswick County voting data \citep{ncsbVoting}. We pre-process the Brunswick County voting data by geocoding the raw addresses to latitude and longitude using the US census bureau geocoding tool and retain all data with a successful geocoding that belong to white and Black voters. For all datasets we normalize all real-valued features to be between $[0,1]$ and convert all categorical features to be real-valued via the one-hot encoding scheme. For datasets that were originally used for supervised learning, we remove the target variable and do not use the sensitive attribute as a feature for the clustering itself. The $k$-medians algorithms require computing a distance matrix between all pairs of points leading to a large memory requirement ($O(n^2)$). To circumvent memory issues for large datasets, we sub-sample all datasets to have at most $10,000$ data points. We note that this memory bottleneck is also present in the baseline $k$-median algorithms \citep{park2009simple} and is independent of our proposed algorithms to the FMRA problem which only requires computing the distance between data points and the candidate centers. We use the same random sub-sample for all algorithms to provide a fair comparison.

We implemented MiniReL in Python with Gurobi 10.0 \citep{gurobi} for solving all IPs. We warm-start MiniReL with the output from the baseline version of Lloyd's algorithm (details and evaluation of this warm-starting is included in Appendix \ref{app:warm_start}). All experiments were run on a computing environment with 16 GB of RAM and 2.7 GHz Quad-Core Intel Core i7 processor. For the following experiments we set $\alpha = 0.51$ to represent majority representation in a cluster. We experiment with different settings of $\alpha$ in Appendix \ref{app:changing_alpha}. Our results show that MR-fairness becomes more difficult to satisfy as the representation threshold $\alpha$ increases, and can be satisfied with standard unfair clustering approaches when $\alpha\le0.3$. We also set  $\ell = 1, u = n$ to provide a fair comparison to Lloyd's algorithm with no cardinality constraints. We provide some additional experiments with balanced clusters (i.e., $l \approx \frac{n}{k}$) in Appendix \ref{app:balanced_clusters}. In the interest of brevity, all experiment results are reported using the cluster statistical parity setting of $\beta$ and we defer experimental results under cluster equality of opportunity to Appendices \ref{app:kmean_extra} and \ref{app:kmed}.

We benchmark MiniReL against Lloyd's algorithm for $k$-means and its associated alternating minimization algorithm for $k$-medians \citep{park2009simple}. For the $k$-means setting we use the implementation available in scikit-learn \citep{scikit-learn} with a $k$-means++ initialization. For the $k$-medians setting we compare MiniReL against the alternating minimization approach for $k$-medians using the implementation available in scikit-learn extra package \citep{scikit-learn-extra} with a $k$-mediods++ initialization. We use Euclidean distance to compute distances between points. For both settings, we run the algorithm with 100 different random seeds. We report the clustering with the lowest clustering cost, which we denote $k$-means/$k$-medians respectively, and the fairest clustering with respect clustering statistical parity ($k$-\textit{means-SP}/$k$-\textit{medians-SP}) and cluster equality of opportunity ($k$\textit{-means-EqOp}/ $k$\textit{-medians-EqOp}). We also benchmark MiniReL against the publicly available implementations of the socially fair $k$-means algorithm of \citet{ghadiri2021socially} in the $k$-means setting, and the diversity aware $k$-medians algorithm of \citet{thejaswi2021diversity}. We do not benchmark against any proportionally fair clustering algorithms \citep[e.g.][]{chierichetti2017fair}, as any algorithm that enforces proportionally fair clustering will be definition violate MR-fairness.

\subsection{Comparing different variants of the MiniReL } \label{sec:algo_design}
\begin{table}[t]
\tabcolsep=0.11cm
\footnotesize
\caption{Impact of algorithmic components on cluster quality, total computation time and average time to solve the assignment problem in the $k$-means setting. Values are normalized to show the percentage change with respect to the base MiniReL algorithm, and averaged over $k \in [2,14]$. Negative values indicate an improvement over the base algorithm.\hfill}.
\resizebox{\linewidth}{!}{\begin{tabular}{llrrrrr}
 data & Metric & MiniReL-TwoStage & MiniReL-Flow & MiniReL-Prefix & MiniReL-PrefixFlow & MiniReL-HeurFlow \\ \toprule
\multirow[c]{4}{*}{\cell{l}{\texttt{heart switzerland}\\ $n=123$, $m=13$ \\ $|{\cal F}|=1$, $|{\cal G}|=2$}} & Avg. Iteration Time & -57.79 & -65.41 & -86.92 & -94.28 & -94.30 \\
 & Fairness Violation & 0.00 & 1.37 & 0.00 & 1.37 & 1.93 \\
 & Objective & -0.80 & -6.19 & -0.80 & -6.19 & -4.02 \\
 & Total Time & -54.06 & 26.16 & -68.05 & -53.50 & -86.19 \\ \cmidrule{1-7}
\multirow[c]{4}{*}{\cell{l}{\texttt{heart va}\\ $n=200$, $m=13$ \\ $|{\cal F}|=1$, $|{\cal G}|=2$}} & Avg. Iteration Time & -68.78 & -75.40 & -89.19 & -94.96 & -94.81 \\
 & Fairness Violation & 0.00 & 1.40 & 0.00 & 1.40 & 1.49 \\
 & Objective & -0.09 & -3.82 & -0.09 & -3.82 & -0.04 \\
 & Total Time & -66.50 & -75.13 & -77.02 & -81.68 & -84.05 \\ \cmidrule{1-7}
\multirow[c]{4}{*}{\cell{l}{\texttt{heart hungarian}\\ $n=294$, $m=13$ \\ $|{\cal F}|=1$, $|{\cal G}|=2$}} & Avg. Iteration Time & -64.86 & -75.07 & -87.55 & -94.87 & -95.04 \\
 & Fairness Violation & 0.00 & 0.53 & 0.00 & 0.53 & 0.63 \\
 & Objective & -0.02 & -2.20 & -0.01 & -2.20 & 0.80 \\
 & Total Time & -63.42 & -74.71 & -79.48 & -85.78 & -90.55 \\ \cmidrule{1-7}
\multirow[c]{4}{*}{\cell{l}{\texttt{heart cleveland}\\ $n=297$, $m=13$ \\ $|{\cal F}|=1$, $|{\cal G}|=2$}} & Avg. Iteration Time & -6.89 & -34.32 & -65.12 & -81.50 & -82.19 \\
 & Fairness Violation & 0.00 & 0.52 & 0.00 & 0.52 & 0.54 \\
 & Objective & -0.01 & -0.71 & -0.01 & -0.71 & -0.28 \\
 & Total Time & -7.84 & 180.33 & -41.26 & 14.38 & -2.59 \\ \cmidrule{1-7}
\multirow[c]{4}{*}{\cell{l}{\texttt{student-mathematics}\\ $n=395$, $m=32$ \\ $|{\cal F}|=1$, $|{\cal G}|=2$}} & Avg. Iteration Time & 65.61 & 56.95 & -18.93 & -6.06 & -5.90 \\
 & Fairness Violation & 0.00 & 0.06 & 0.00 & 0.06 & 0.05 \\
 & Objective & -0.00 & 0.01 & -0.00 & 0.01 & -0.01 \\
 & Total Time & 74.51 & 16.88 & 54.24 & 59.34 & 31.69 \\ \cmidrule{1-7}
\multirow[c]{4}{*}{\cell{l}{\texttt{house-votes}\\ $n=435$, $m=16$ \\ $|{\cal F}|=1$, $|{\cal G}|=2$}} & Avg. Iteration Time & 28.51 & 14.89 & -52.88 & -58.98 & -55.80 \\
 & Fairness Violation & 0.00 & 0.19 & 0.00 & 0.19 & 0.16 \\
 & Objective & -0.01 & 0.21 & -0.00 & 0.21 & 1.44 \\
 & Total Time & 24.11 & 1.10 & -18.15 & -21.55 & -7.49 \\ \cmidrule{1-7}
\multirow[c]{4}{*}{\cell{l}{\texttt{student-por}\\ $n=649$, $m=32$ \\ $|{\cal F}|=1$, $|{\cal G}|=2$}} & Avg. Iteration Time & 34.73 & 13.64 & -55.54 & -38.21 & -36.02 \\
 & Fairness Violation & 0.00 & 0.27 & 0.00 & 0.27 & 0.23 \\
 & Objective & 0.00 & -0.05 & 0.03 & -0.05 & -0.04 \\
 & Total Time & 60.12 & 45.79 & -17.95 & 22.52 & 16.20 \\ \cmidrule{1-7}
\multirow[c]{4}{*}{\cell{l}{\texttt{hcv}\\ $n=1385$, $m=28$ \\ $|{\cal F}|=1$, $|{\cal G}|=2$}} & Avg. Iteration Time & 29.18 & 19.88 & -66.84 & -58.64 & -49.00 \\
 & Fairness Violation & 0.00 & 0.07 & 0.00 & 0.07 & 0.08 \\
 & Objective & 0.00 & -0.00 & 0.00 & -0.00 & -0.00 \\
 & Total Time & 22.69 & 4.87 & 7.86 & 19.27 & -15.56 \\ \cmidrule{1-7}
\multirow[c]{4}{*}{\cell{l}{\texttt{abalone}\\ $n=4177$, $m=8$ \\ $|{\cal F}|=1$, $|{\cal G}|=3$}} & Avg. Iteration Time & -64.64 & -74.10 & -95.80 & -96.82 & -96.83 \\
 & Fairness Violation & 0.00 & 0.12 & 0.00 & 0.12 & 0.11 \\
 & Objective & 0.00 & -0.04 & 0.00 & -0.04 & 0.10 \\
 & Total Time & -61.33 & -70.61 & -71.53 & -64.04 & -90.74 \\ \cmidrule{1-7} \cmidrule{1-7}
\multirow[c]{4}{*}{\cell{l}{\texttt{default}\\ $n=30000$, $m=20$ \\ $|{\cal F}|=1$, $|{\cal G}|=2$}} & Avg. Iteration Time & -44.05 & -48.27 & -93.27 & -93.88 & -93.91 \\
 & Fairness Violation & 0.00 & 0.06 & 0.00 & 0.05 & 0.04 \\
 & Objective & 0.00 & -0.01 & 0.00 & -0.01 & 0.09 \\
 & Total Time & -41.40 & -45.48 & -42.82 & -40.11 & -84.44 \\ \cmidrule{1-7}
\multirow[c]{4}{*}{\cell{l}{\texttt{adult}\\ $n=32561$, $m=11$ \\ $|{\cal F}|=2$, $|{\cal G}|=5$}} & Avg. Iteration Time & -39.31 & -40.61 & -95.64 & -96.21 & -96.10 \\
 & Fairness Violation & 0.00 & 0.06 & 0.00 & 0.06 & 0.07 \\
 & Objective & -0.14 & -0.12 & -0.23 & -0.22 & 1.25 \\
 & Total Time & -51.38 & -55.60 & -71.84 & -72.10 & -90.47 \\ \cmidrule{1-7}
\multirow[c]{4}{*}{\cell{l}{\texttt{voting}\\ $n=49190$, $m=3$ \\ $|{\cal F}|=1$, $|{\cal G}|=2$}} & Avg. Iteration Time & -37.45 & -43.88 & -94.12 & -95.10 & -95.07 \\
 & Fairness Violation & 0.00 & 0.02 & 0.00 & 0.02 & 0.02 \\
 & Objective & -4.48 & -4.51 & -5.13 & -5.15 & 1.67 \\
 & Total Time & -39.11 & -45.14 & -79.44 & -80.61 & -93.02 \\ \cmidrule{1-7}
\multirow[c]{4}{*}{\cell{l}{\texttt{diabetes}\\ $n=95309$, $m=2$ \\ $|{\cal F}|=2$, $|{\cal G}|=4$}} & Avg. Iteration Time & -87.70 & -96.97 & -93.16 & -98.14 & -98.14 \\
 & Fairness Violation & 0.00 & 0.01 & 0.00 & 0.01 & 0.01 \\
 & Objective & -7.03 & -7.04 & -7.03 & -7.04 & -8.02 \\
 & Total Time & -25.72 & -80.27 & -55.06 & -83.61 & -88.41 \\
\end{tabular}
}
\label{tab:big_table}
\end{table}

We evaluate the following 5 different variants of the MiniReL algorithm to show the impact of different algorithmic components on its performance:
\begin{itemize}

    \item MiniReL-TwoStage: Decomposes the FMRA by sequentially solving the GRP, by relaxing the $z$ variables and solving the relaxed problem via IP, and then solving APFR, using IP.
    \item MiniReL-Flow: Uses the two-stage decomposition but solves the APFR via the network flow approach.
    \item MiniReL-Prefix: Uses the two-stage decomposition solving the GRP, by relaxing the $z$ variables and solving the problem via IP once, and solves the second stage APFR, using IP.
    \item MiniReL-PrefixFlow: Uses the two-stage decomposition with pre-fixing (i.e., solves the \fsabbrv via IP once) and solves the second stage using the network flow approach introduced in Section \ref{sec:flow}.
    \item MiniReL-HeurFlow: Uses the two-stage decomposition with pre-fixing with the heuristic introduced in Section \ref{sec:prefix} and solves the second stage using the network flow approach introduced in Section \ref{sec:flow}.
\end{itemize}

For the sake of brevity we report results for this ablation study in the $k$-means setting, but the results were similar in the $k$-medians setting. Table \ref{tab:big_table} shows the performance of five algorithmic variants with respect to the average iteration time (i.e., the time to do one iteration of assignment and center computation), the normalized additive fairness violation (i.e., the total additive fairness violation normalized by the size of the dataset), cluster quality (Objective), and total computation time of the algorithm. All values presented, with the exception of the additive fairness violation, are normalized to show the percentage change with respect to the base MiniReL algorithm (i.e., solving the full FMRA to optimality in each iteration), and averaged over $k \in [2,14]$. Negative values indicate an improvement over the baseline algorithm. All algorithms were given a 10 minute time limit, and the results display the quality of the best clustering found within the time limit. 

Overall, the results validate key algorithmic components of MiniReL. For large datasets, the two-stage decomposition (MiniReL-TwoStage) approach leads to a decrease in computation time over the baseline algorithm with as much as a 50\% reduction in computation time on the Adult dataset. For smaller datasets the two-stage decomposition has a mixed effect, increasing the computation time in some instances while reducing it in others. For smaller datasets, the baseline MiniReL approach runs quickly and the two-stage decomposition introduces additional overhead (i.e., setting up two MILP problems) that outweights its benefits. Solving the second stage problem via network flow (MiniReL-Flow) also gives a small additional reduction in computation time over the two-stage decomposition alone, especially on a per-iteration basis. Both approaches have no large impact on the objective (i.e., clustering cost) of the final solutions, and the addition of the network flow algorithm only adds very modest additive fairness violations (at most 1.4\% for small datasets and 0.06\% for larger datasets). Adding pre-fixing (Minirel-Prefix) yields an additional reduction in computation time in most of the datasets and dramatically reduces the per iteration computation time. Solving the \fsabbrv problem via the heuristic led to a large additional reduction in computation time (e.g., an additional $20\%$ in the adult dataset) but comes with a slight reduction in the quality of the clusters found (at worst 2\% increase in the objective). We also validate these results on synthetic data in Appendix \ref{app:synthetic}.

Based on these results we focus on the two most promising variants of our algorithm: MiniReL-Prefix, which we denote with MiniReL-IP, and MiniReL-PrefixFlow, which we denote with MiniReL-Flow. We evaluate both variants as they include the most successful algorithmic components, pre-fixing and heuristics for the \fsabbrvx. While MiniReL-IP satisfies MR-fairness exactly,  MiniReL-Flow incorporates the faster bi-criteria network approximation for the APFR and might have small additive violations of the fairness constraints.

\subsection{k-means results} \label{sec:kmeans_exp}

\begin{figure}[!thb]
    \centering
    \begin{subfigure}
      \centering
      \includegraphics[width=0.98\textwidth]{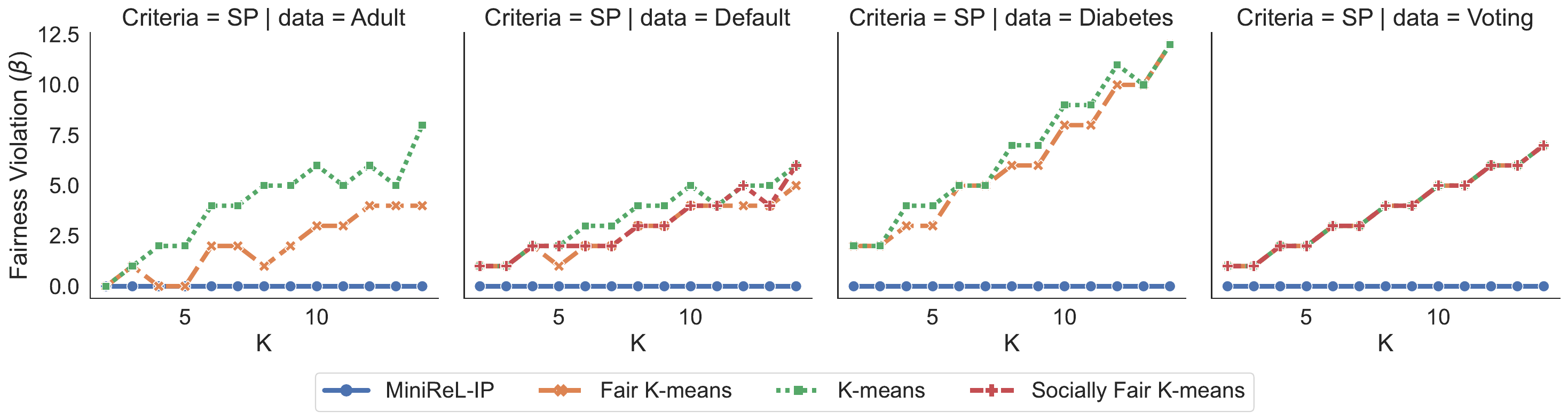}
    \end{subfigure}
    \begin{subfigure}
      \centering
      \includegraphics[width=0.98\textwidth]{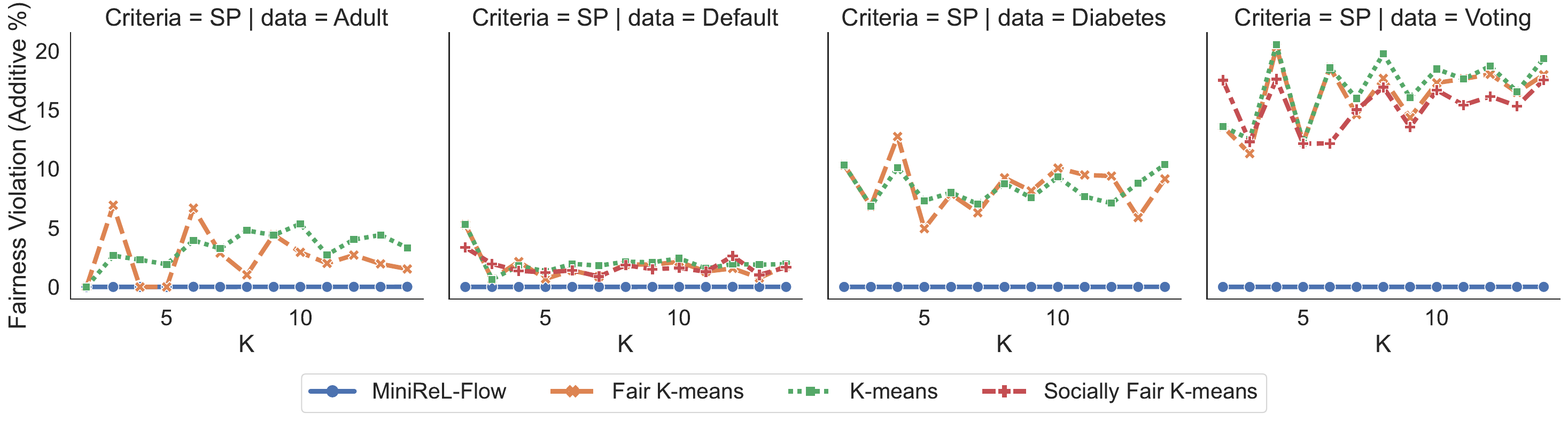}
    \end{subfigure}
  \caption {\label{fig:fairness_joint} \vspace{-0.1cm}(Top) Maximum fairness violation for MiniReL-IP and baseline algorithms in $k$-means setting under statistical parity. (Bottom) Normalized Additive Fairness Violation (additive fairness violation divided by dataset size) for MiniReL-Flow and baseline algorithms under statistical parity.
  }
\end{figure}
We compare MiniReL-Flow and MiniReL-IP in the $k$-means setting with standard Lloyd's algorithm for $k$-means and the socially fair k-means algorithm \citep{ghadiri2021socially}. We show results for the four largest datasets from our extended benchmark suite --- Adult, Default, Diabetes, and Voting. Note that the socially fair k-means algorithm is only designed to work with a single sensitive feature and thus we only report results for the Default and Diabetes datasets. We evaluate all algorithms with respect to cluster statistical parity (i.e., every group must be $\alpha$-represented in the same number of clusters), and report results for cluster equality of opportunity in Appendix \ref{app:kmean_extra}.

Figure \ref{fig:fairness_joint} (top row) shows the maximum fairness violation (i.e., $\max_{g \in {\mathcal G} }\max(\beta_g - \Lambda (C, X_g, \alpha), 0)$) for the baseline algorithms and MiniReL-Exact.  Figure \ref{fig:fairness_joint} (bottom row) shows the total additive fairness violation normalized by the size of the dataset for the baseline algorithms and MiniReL. 
Across all four datasets we can see that $k$-means can lead to outcomes that violate MR-fairness constraints significantly, and that selecting the fairest clustering has only a marginal improvement. This is most striking in the default dataset where there is as much as an 11 cluster gap (6 cluster fairness violation) between the two groups despite having similar proportions in the dataset. Similarly, enforcing a related notion of fairness, Social Fairness, has no noticeable improvement on the fairness of the clustering under MR-Fairness, empirically validating the theoretical results from Section \ref{sec:fair_metric_comp}. In contrast, MiniReL-IP is able to generate fair clusters by design. This result also applies in the additive fairness violation setting, where the baseline algorithms generate large constraint violations (as much as 20\% of the dataset size for the voting dataset) whereas MiniReL-Flow has normalized additive fairness violations close to 0\%. 

Table \ref{table:comp_time} shows the average computation time in seconds for Lloyd's algorithm and MiniReL-Flow under cluster statistical parity. We display results for the cluster statistical parity as it was the more computationally demanding setting requiring more MiniReL-Flow iterations than cluster equality of opportunity. As expected, the harder assignment problem in MiniReL-Flow leads to higher overall computation times when compared to the standard Lloyd's algorithm. However, MiniReL-Flow is still able to solve large problems in under 200 seconds demonstrating that the approach is of practical use.

\begin{table}[t]
\tabcolsep=0.11cm
\centering
\footnotesize
\caption{Average computation time (standard deviation) in seconds for MiniReL-Flow under Statistical Parity in the $k$-means setting. Final row includes average computation time over all settings of K for each dataset. \label{table:comp_time}}.
\begin{tabular}{l|ll|ll|ll|ll}
\toprule
data & \multicolumn{2}{c|}{Adult} & \multicolumn{2}{c|}{Default} & \multicolumn{2}{c|}{Diabetes} & \multicolumn{2}{c}{Voting} \\
k & K-means & MiniReL-Flow & K-means & MiniReL-Flow & K-means & MiniReL-Flow & K-means & MiniReL-Flow \\
\midrule
2 & 1.91 & 13.38 & 0.65 & 17.64 & 1.34 & 19.54 & 0.48 & 41.71 \\
3 & 2.6 & 22.91 & 0.82 & 15.0 & 1.29 & 19.95 & 0.56 & 37.92 \\
4 & 3.0 & 37.57 & 0.79 & 16.33 & 1.14 & 23.4 & 1.3 & 79.81 \\
5 & 3.59 & 41.06 & 1.25 & 17.5 & 1.44 & 35.19 & 0.84 & 58.37 \\
6 & 5.16 & 24.58 & 1.42 & 19.15 & 1.0 & 38.52 & 0.92 & 127.19 \\
7 & 3.66 & 26.49 & 2.36 & 25.1 & 1.66 & 46.56 & 1.07 & 67.27 \\
8 & 5.7 & 27.93 & 1.5 & 23.43 & 1.7 & 47.27 & 1.48 & 113.66 \\
9 & 6.07 & 60.11 & 1.53 & 24.11 & 1.72 & 62.86 & 1.53 & 69.07 \\
10 & 6.33 & 71.18 & 1.88 & 24.69 & 2.31 & 69.79 & 2.76 & 81.41 \\
11 & 3.79 & 74.84 & 2.01 & 26.45 & 1.69 & 50.75 & 1.9 & 198.02 \\
12 & 5.24 & 74.52 & 2.66 & 28.72 & 3.01 & 58.05 & 2.93 & 108.86 \\
13 & 4.08 & 60.11 & 3.21 & 26.95 & 2.12 & 51.53 & 3.23 & 186.94 \\
14 & 2.91 & 64.61 & 8.42 & 28.28 & 1.6 & 79.28 & 3.21 & 218.53 \\ \midrule
\textbf{mean} & \textbf{4.16} & \textbf{46.10} & \textbf{2.19} & \textbf{22.57} & \textbf{1.69} & \textbf{46.36} & \textbf{1.71} & \textbf{106.83} \\
\bottomrule
\end{tabular}
\end{table}

\begin{figure}[b]
    \centering
    \begin{subfigure}
      \centering
      \includegraphics[width=\textwidth]{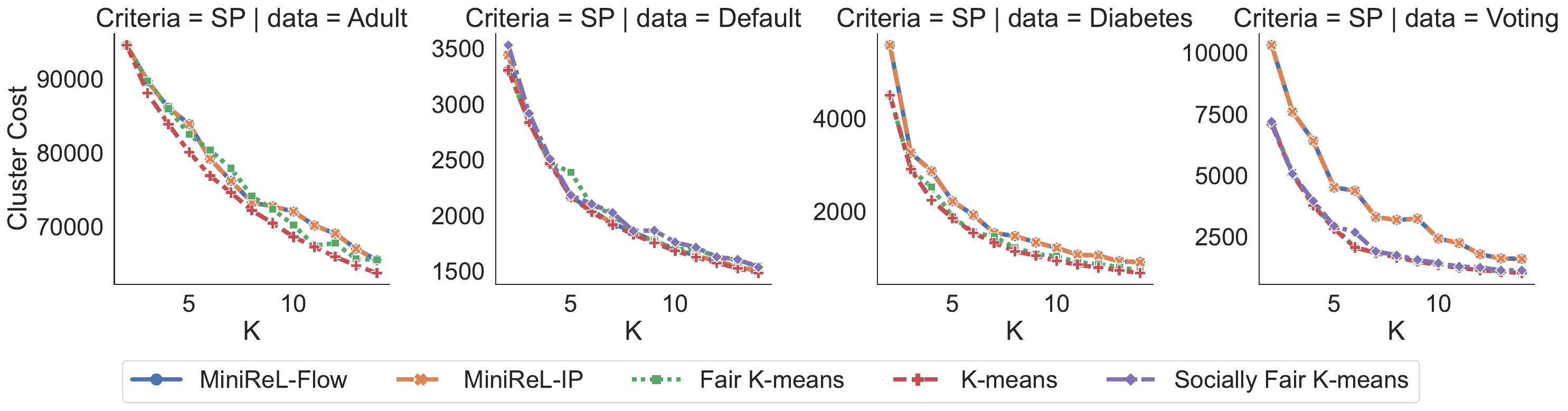}
    \end{subfigure}
    \vspace{-0.5cm}
  \caption {\label{fig:cluster_cost_joint} $k$-means clustering cost of fair clustering algorithms under statistical parity.
  }
\end{figure}

 One remaining question is whether fairness comes at the expense of cluster quality. Figure \ref{fig:cluster_cost_joint} shows the $k$-means clustering cost for the baseline algorithms and MiniReL. Although there is a small increase in the cost when using MiniReL the overall cost closely matches that of the standard $k$-means algorithm in three out of four datasets, with the sole exception being the voting dataset, showing that we can gain fairness at practically no additional increase to clustering cost. 

\subsection{k-medians results} \label{sec:kmedians_exp}
In this section we benchmark MiniReL-Flow and MiniReL-IP against the Lloyd-style $k$-medians algorithm \citep{park2009simple} and the diversity-awware $k$-medians algorithm \citep{thejaswi2021diversity}. We present results for the algorithm under statistical parity and include the extended results in Appendix \ref{app:kmed}. Figure \ref{fig:kmed_fairness_sp} shows the maximum fairness violation and normalized additive fairness violation for MiniReL and the baseline algorithms. Similar to the $k$-means setting, the $k$-medians algorithm without fairness constraints can lead to unfair outcomes across all datasets, whereas MiniReL-IP constructs fair clusters by design. Figure \ref{fig:kmed_cluster_cost_sp} shows the cluster cost for both algorithms in the same setting. Adding fairness in these settings again leads to only moderate increase in cluster cost, with the only notable increase coming on the voting dataset under statistical parity. 

Table \ref{tab:kmed_runtime} shows the runtime of both algorithms under statistical parity. MiniReL-Flow still leads to a modest increase in runtime compared to the baseline algorithm but is able to solve all datasets (10,000 data points after sub-sampling) in under 40 seconds in most instances.

\begin{figure}[!t]
    \centering
    \begin{subfigure}
      \centering
      \includegraphics[width=0.9\textwidth]{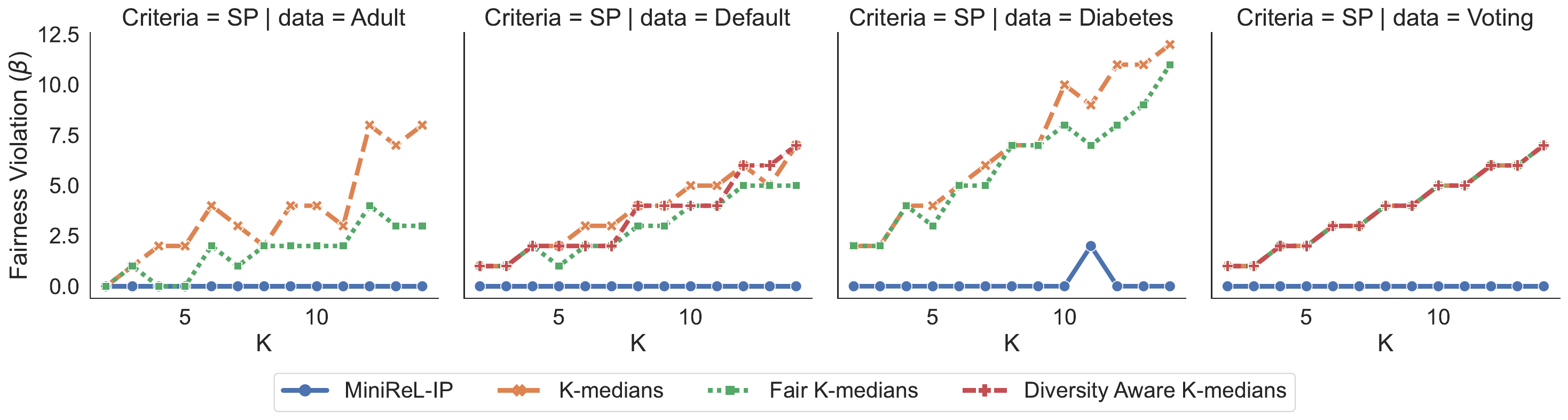}
    \end{subfigure}
  \caption {\label{fig:kmed_fairness_sp} Maximum fairness violation for Lloyd's algorithm and MiniReL in $k$-medians setting with IP assignment under Statistical Parity.
  }
\end{figure}

\begin{figure}[!t]
    \centering
    \begin{subfigure}
      \centering
      \includegraphics[width=\textwidth]{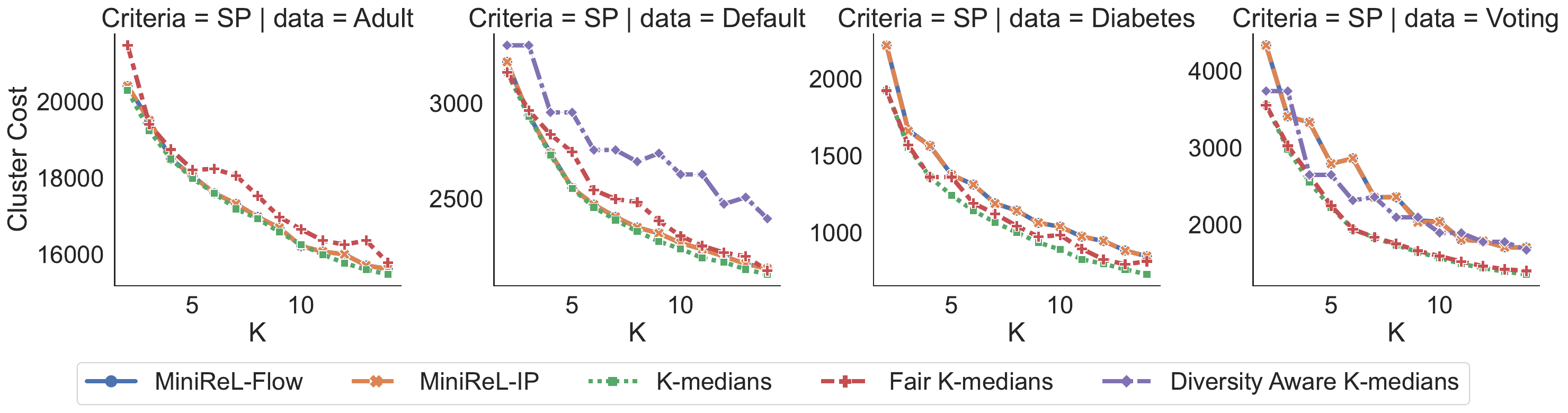}
    \end{subfigure}
    \vspace{-0.5cm}
  \caption {\label{fig:kmed_cluster_cost_sp} $k$-medians clustering cost of Lloyd's algorithm and MiniReL under statistical parity.
  }
\end{figure}

\begin{table}
\tabcolsep=0.11cm
\footnotesize 
\centering
\caption{Average computation time in seconds for MiniReL-Flow under Statistical Parity in $k$-medians setting. Final row includes average computation time over all settings of K for each dataset. Note all datasets are sub-sampled to 10,000 points \label{tab:kmed_runtime}}.
\begin{tabular}{l|rr|rr|rr|rr}
\toprule
data & \multicolumn{2}{c|}{Adult} & \multicolumn{2}{c|}{Default} & \multicolumn{2}{c|}{Diabetes} & \multicolumn{2}{c}{Voting} \\
k & K-medians & MiniReL-Flow & K-medians & MiniReL-Flow & K-medians & MiniReL-Flow & K-medians & MiniReL-Flow \\
\midrule
2 & 75.43 & 16.02 & 103.42 & 17.24 & 37.71 & 4.66 & 127.95 & 20.63 \\
3 & 41.39 & 14.31 & 67.76 & 10.53 & 17.96 & 5.5 & 72.85 & 21.37 \\
4 & 33.23 & 10.15 & 86.82 & 16.01 & 13.74 & 4.41 & 42.33 & 17.99 \\
5 & 51.81 & 10.58 & 60.69 & 15.89 & 9.95 & 4.12 & 49.65 & 14.69 \\
6 & 33.61 & 11.67 & 32.96 & 11.4 & 11.03 & 2.84 & 41.95 & 23.76 \\
7 & 29.02 & 10.57 & 44.89 & 11.7 & 9.88 & 3.11 & 56.5 & 22.65 \\
8 & 24.96 & 11.93 & 42.21 & 11.89 & 5.34 & 3.95 & 59.96 & 23.3 \\
9 & 27.93 & 11.91 & 34.98 & 10.92 & 6.88 & 4.99 & 31.68 & 22.38 \\
10 & 26.33 & 20.0 & 37.99 & 12.81 & 7.29 & 4.2 & 27.3 & 23.94 \\
11 & 20.02 & 12.49 & 27.56 & 12.48 & 4.53 & 5.7 & 31.32 & 35.91 \\
12 & 15.1 & 12.54 & 29.27 & 13.2 & 4.52 & 5.04 & 26.5 & 32.77 \\
13 & 14.59 & 13.59 & 34.06 & 17.61 & 4.32 & 3.85 & 36.16 & 28.66 \\
14 & 19.82 & 14.35 & 43.34 & 14.3 & 4.5 & 6.48 & 29.63 & 39.28 \\ \midrule
\textbf{mean} & \textbf{31.79} & \textbf{13.08} & \textbf{49.69} & \textbf{13.54} & \textbf{10.59} & \textbf{4.53} & \textbf{48.75} & \textbf{25.18} \\
\bottomrule
\end{tabular}
\end{table}

\section{Conclusion}
In this paper we introduce a novel definition of group fairness for clustering that requires each group to have a minimum level of representation in a specified number of clusters. This definition is a natural fit for a number of real world examples, such as voting and entertainment segmentation. To create fair clusters we introduce a modified version of Lloyd's algorithm called MiniReL that solves an assignment problem in each iteration via integer programming. While the underlying optimization problem is computationally challenging, we show that MiniReL can be implemented efficiently through a combination of decomposition, network flow formulations, and heuristic pre-fixing strategies. Extensive experiments across real and synthetic datasets demonstrate that MiniReL consistently produces fair clusterings while maintaining cluster quality and computational tractability.

\newpage
\bibliographystyle{or_style/informs2014}
\bibliography{references}

\newpage
\begin{APPENDIX}{}
\section{Table of Notation} \label{app:notation}

\begin{table*}[h!]

\begin{tabular}{l|l}
\textbf{Symbol} & \textbf{Definition} \\ \hline
\addlinespace[1mm]\multicolumn{2}{c}{Clustering Notation}\\
\hline
$\cal{X} = \{x^i \in \mathbb{R}^m\}_{i=1}^n$  & Dataset of $n$ $m$-dimensional data points\\
$K$ & The number of clusters (given as input to the problem) \\
${\cal K} = \{1, \dots, K\}$ & The index set of clusters \\
${\cal C} = \{C_1, \dots, C_k\}$ & The set of $K$ clusters \\
${\cal L}_i$ & The set of allowable centers for the $i$-th cluster \\
$c_i \in {\cal L}_i$ & The center of a cluster $i$ \\
$D(x^i, c)$ & Cluster cost for assigning data point $i$ to center $c$ \\
$C_i \subseteq {\cal X} $ & The set of data points in cluster i \\
${\cal F}$ & The set of sensitive features (e.g., Gender, Race) \\
$u$ & Upper bound on the cardinality of each cluster\\
$l$ & Lower bound on the cardinality of each cluster\\
${\cal G}_f \subseteq {\cal F}$ & The set of possible values for a sensitive feature $f$ \\
${\cal G} = \bigcup_{f \in {\cal F}} {\cal G}_f$ & The set of possible groups (each group corresponds to one value of a sensitive feature) \\
$X_g$ & The set of data points belonging to a group $g$ \\
\hline
\addlinespace[1mm]\multicolumn{2}{c}{Fairness Notation}\\
\hline
$\alpha$ & The representation threshold for MR-fair clustering \\
$\mathbf{\beta}$ & The required number of clusters with $\alpha$-representation for each group \\
$\Lambda({\cal C}, X_g, \alpha) $ & The number of clusters that have $\alpha$-representation for group $g$ in a clustering ${\cal C}$ \\
\hline
\addlinespace[1mm]\multicolumn{2}{c}{IP Formulation Notation}\\
\hline
$z_{ik} \in \{0,1\}$ & Binary variable indicating whether data point $i$ assigned to cluster $k$ \\
$y_{gk} \in \{0,1\}$ & Binary variable indicating whether group $g$ needs to be $\alpha$-represented in cluster $k$\\
$m_{gk}$ & The myopic cost of $\alpha$-representing group $g$ in cluster $k$\\
${\cal T}$ & The set of \textit{types}, a combination of sensitive features (e.g., Male and Black)\\

\hline
\addlinespace[1mm]\multicolumn{2}{c}{Acronyms}\\
\hline
FMRA & The fair minimum-representation assignment problem ( \eqref{obj:mip}-\eqref{const:binary_y} with fixed centers) \\
GRP & The group representation problem ( \eqref{obj:mip}-\eqref{const:binary_y} with fixed centers, $z_{ik} \in [0,1]$) \\
APFR & The assignment problem under fixed representation ( \eqref{obj:mip}-\eqref{const:binary_y} with fixed centers, $y_{gk}$) \\
Flow-APFR & The APFR problem modeled as the network flow problem described in Section \ref{sec:flow}\\
\hline
\end{tabular}
\caption{List of common notation}
\label{table:notation}
\end{table*}

\section{Proof of Theorem \ref{thm:fair_assign_np}} \label{app:fair_assign_np}
\newcommand{\calu}{{\cal U}}\newcommand{\calw}{{\cal W}}

\proof{Proof}
The reduction is from the exact cover by 3 sets  (X3C) problem, 
with $3q$ elements ${\cal U} = \{u_1, \dots, u_{3q}\}$ and $t = q + r$ subsets ${\cal W} = \{W_1, \dots, W_t\}$ where each subset $W_i \subseteq {\cal U}$ and $|W_i| = 3$. 
Recall that X3C is one of Karp's 21 NP-Complete problems and its goal  is to select a collection of subsets ${\cal W}^* \subseteq {\cal W}$ such that all elements of ${\cal U}$ are covered exactly once.

Given an instance of the X3C problem, we first construct an undirected bipartite graph $G = ({\calu\cup\calw}, E)$ which we then convert into a FMRA instance following a similar construction as \citep{esmaeili2022fair}. 
Graph $G$ has a node for each $u_i\in\calu$ and a node $w_j$ for each  $W_j\in\calw$. There is an edge $\{u_i,w_j\}$ between nodes
$u_i$ and ${w_j} $ provided that $u_i \in W_j$. Note that with slight abuse of notation, we use $\calw$ both for the subsets in the X3C instance and the nodes corresponding to them in the graph $G$.

Using graph $G$, we construct an instance of the FMRA problem with two groups by interpreting each vertex in the graph as a data point where points corresponding to vertices in ${\cal U}$ and ${\cal W}$ belong to separate groups.
We set  $K = t$ (i.e., one cluster for each subset), $\alpha = 0.75$,  $\beta_W = t - q$  and $\beta_U = q$ (minimum number of $\alpha$-represented clusters for each group). 
We place the (fixed) $K$ centers at the data points associated with $\cal W$. 
The distances between pairs of data points are set to be the length of the shortest path distances between the corresponding vertices in the graph.

We now argue that the optimal solution to the FMRA instance has a cost of $3q$ if and only if there exists a feasible solution to the X3C instance. 
First note that any feasible solution  to the FMRA instance has a cost of at least $3q$ as each point $u_i\in\calu$ incurs a cost of at least 1 and there are $3q$ such points. Next, note that given a feasible solution $\calw^*\subseteq\calw$  to the X3C instance, we can construct a solution to the corresponding FMRA instance with cost $3q$ by creating one cluster for each $W_j\in\calw^*$ in the solution.
These $q$ clusters contain $w_j$ together with the 3 elements $u_i\in W_j$. The remaining $t-q$ clusters contain a single point $w_j$, one for for each $W_j\in\calw\setminus\calw^*$. Note that these clusters satisfy the fairness constraints. Furthermore, letting the  center of each cluster to be the $w_j$ that it contains gives a cost of $3q$.

Finally we argue that if the FMRA instance has a cost of $3q$, then the X3C instance is feasible. If the solution to the FMRA has cost $3q$, then $(i)$ each point $u_i\in\calu$ must be assigned to a center $w_j\in\calw$ such that $u_i\in W_j$, and, $(ii)$ each point $w_j\in\calw$ must be assigned to the center $w_j$ and consequently, each cluster must contain exactly one $w_j\in\calw$.
In addition, the fairness constraint for the $\calw$ group requires that at least $t-q$ clusters must have 75\% of its points belonging to the $\calw$ group. As there is exactly one $w_j\in\calw$ in each cluster, this can only happen if at least $t-q$ clusters have no points from $\calu$ group in them.
Consequently all $3q$ points form the $\calu$ group are assigned to at most $q$ clusters. As all $u_i\in\calu$ must be assigned to a center $w_j\in\calw$ such that $u_i\in W_j$, and $|W_j|=3$, the solution to the FMRA instance must have exactly $t-q$ clusters with no points from $\calu$ group and $q$ clusters with one $w_j$ and three $u_i\in W_j$, which gives a feasible solution to the X3C instance. \qed
\endproof
	
\section{Proof of Theorem \ref{thm:finite_convergence}} \label{app:pf_convergence}
\proof{Proof.}
We start by noting that given an assignment of data points to clusters, the optimal cluster centers for this assignment are precisely the ones computed as in Step \ref{alg:center_mean} (\ref{alg:center_med}).
Consequently, if \textit{improved\_cost = current\_cost} in Step \ref{alg:until}, then the cluster centers used in Step \ref{alg:assign} and the ones computed in Step \ref{alg:center_mean} or \ref{alg:center_med} (in the last iteration) must be identical. This follows from the uniqueness of the optimal centers in the $k$-means setting and the deterministic tie-breaking in the $k$-medians setting. Therefore, if Algorithm 1 terminates, then the current assignment of the data points to the current cluster centers is optimal and the cost cannot improve by changing their assignment due to Step \ref{alg:assign}.
In addition, Steps \ref{alg:center_mean} and \ref{alg:center_med} guarantees that perturbing cluster centers cannot not improve the cost for the current assignment.

It remains to show that the algorithm will terminate after a finite number of iterations. 
Similar to the proof for Lloyd's algorithm, we leverage the fact that there are only a finite number of partitions of the data points. By construction, the objective value decreases in each iteration of the algorithm, and thus can never cycle through any partition multiple times as for a given partition we use the optimal cluster centers when computing \textit{improved\_cost} in Step \ref{alg:improved_cost}. Thus the algorithm can visit each partition at most once and thus must terminate in finite time.
\qed

\endproof
		
\section{Proof of Theorem \ref{thm:second_stage_np}}
\label{app:second_stage_np}

\proof{Proof} Given an instance of a 3-SAT problem, one of Karp's 21 NP-Complete problems, we construct an instance of the \ssabbrvx~problem as follows. We start with a 3-SAT problem with $n$ variables and $m$ clauses $K_1, \dots, K_m$. Each clause $K_i$ takes the form $v_{i1} \vee v_{i2} \vee v_{i3}$ where $v_{ij}$ is either one of the original variables or its negation. 

To construct the instance of the \ssabbrvx~problem, start by creating two new data points $x_{v_i}, x_{\bar{v}_i}$, corresponding to each original variable $v_i$ and its negation respectively. We construct two clusters $C_1$ and $C_2$ that each variable can be assigned to. Let ${\cal W}$ represent the set of allowable group-cluster assignments (i.e., $(g, k) \in {\cal W}$ if $y_{gk} = 1$). For each original variable $v_i$ we create one group $g_i = \{x_{v_i}, x_{\bar{v}_i}\}$ that can be $\alpha$-represented in either cluster (i.e. $(g_i, 1), (g_i,2) \in \mathcal{W}$). For each group $g_i$ we set $\beta_{g_i} = 2$ - ensuring that both clusters must be $\alpha$-represented by the group. We also create one group for each clause $K_i$ corresponding to its three conditions $g_{K_i} = \{x_{v_{i1}}, x_{v_{i2}}, x_{v_{i3}}\}$ that must be $\alpha$-represented in cluster $1$ (i.e. $(g_{K_i},1) \in \mathcal{W}, (g_{K_i},2) \notin \mathcal{W}$). For these groups we set $\beta_{K_i} = 1$. Finally, we set $\alpha = \frac{1}{2n}$ - this ensures that any assignment of a group's data point to a cluster will satisfy the $\alpha$-representation constraint (as there are $2n$ data points and thus at most $2n$ data points in a cluster). We also add a cardinality lower bound on both clusters of 1. Clearly the above scheme can be set-up in polynomial time.

We now claim that a feasible solution to the aforementioned \ssabbrvx~problem corresponds to a solution to the original 3-SAT instance. We start by taking the variable settings by looking at $C_1$. We start by claiming that for each variable exactly one of $x_{v_i}, x_{\bar{v}_i}$ are included in $C_1$. Suppose this weren't true, either both variables were included in $C_1$ or $C_2$. However, whichever cluster has neither of the variables would not be $\alpha$-represented by group $g_i$, thereby contradicting the fairness constraints. Since $C_1$ contains either $x_{v_i}$ or $x_{\bar{v}_i}$ we set $v_i = T$ if $x_{v_i}$ is included and $v_i = F$ otherwise. We now claim that such a setting of the variables satisfies all the clauses. Assume it did not, then there exists a clause $K_i$ such that none of $x_{v_{i1}}, x_{v_{i2}}, x_{v_{i3}}$ are included in $C_1$. However, this violates the $\alpha$-representation constraint for $g_{K_i}$ providing a contradiction to the feasibility of the fair assignment problem. Note that since the feasibility problem does not consider the objective or the location of the centers the proof applies for both the $k$-means and $k$-medians setting. \qed
		
\section{Experiments with different objectives for \fsabbrvx~heuristic} 
\label{app:prefix_obj}

To test the efficacy of the myopic cost for pre-fixing we perform a computational study. We experiment with three different choices of cost function:
\begin{itemize}
    \item \textit{Proportion}: Set the cost to the proportion of the cluster that needs to be changed for $g$ to be $\alpha$-represented in cluster $k$: $$c_{gk} = \max(\alpha - p_{gk},0)$$ where $p_{gk}$ is the current proportion of cluster $k$ that belongs to group $g$.
    \item \textit{Weighted Proportion}: Set the cost to the proportion weighted by the size of the cluster: $$ c_{gk} = |C_k|\max(\alpha - p_{gk},0).$$
    \item \textit{Local Cost}: The myopic (i.e. one-shot) increase in cost of moving additional data points to cluster $k$ to meet the fairness constraints. Let $q$ be the additional number of points from group $g$ needed for $g$ to be $\alpha$-represented in this cluster (i.e. $q + p_{gk}|C_k| \geq \alpha(q + |C_k|)$). Let $c_{(x)} = \argmin_{ c \in \{c_1, \dots, c_K\}} \{D(x, c)\}$ be the closest center for $x \in {\cal X}$. The myopic cost $c_{gk}$ associated with $\alpha$-representing group $g$ in cluster $k$ then becomes: 
    $$
     c_{gk} = \min_{X \subset X_g \setminus C_K: |X| = q} \sum_{x \in X} \big(D(x,c_k) - D(x, c_{(x)})\big)
    $$
\end{itemize}

Figure \ref{fig:prefix_time} shows the impact of different choices of the objective in IP model (8)-(11) on run-time for MiniReL on three different datasets. For these experiments we focused on the $k$-means setting (though $k$-medians showed a similar result) and ran 10 trials with different random initial seeds, and warm-start the algorithm with the standard Lloyd's algorithm. We benchmark using the IP model to perform pre-fixing with naively pre-fixing the group cluster assignments (i.e. random assignment), and compare their performance on three datasets from the UCI machine learning repository \cite{dua2017uci}: Iris, Adult, and Default. In the small 150 data point iris dataset, the pre-fixing scheme has little impact on the total run-time of the algorithm as the overhead of running the IP model outweighs any time savings from a reduced number of iterations. However, for larger datasets (i.e. adult and default which both have over 30K data points), using the IP model to perform pre-fixing outperforms the naive approach - leading to as large as a 3x speed-up. However, there is a relatively small difference in performance between the three choices for the objective, with the local cost objective reducing the speed by approximately $0.1\%$ compared to the other two when averaged across all three datasets. 

\begin{figure}[b]
    \centering
    \begin{subfigure}
      \centering
      \includegraphics[width=\textwidth]{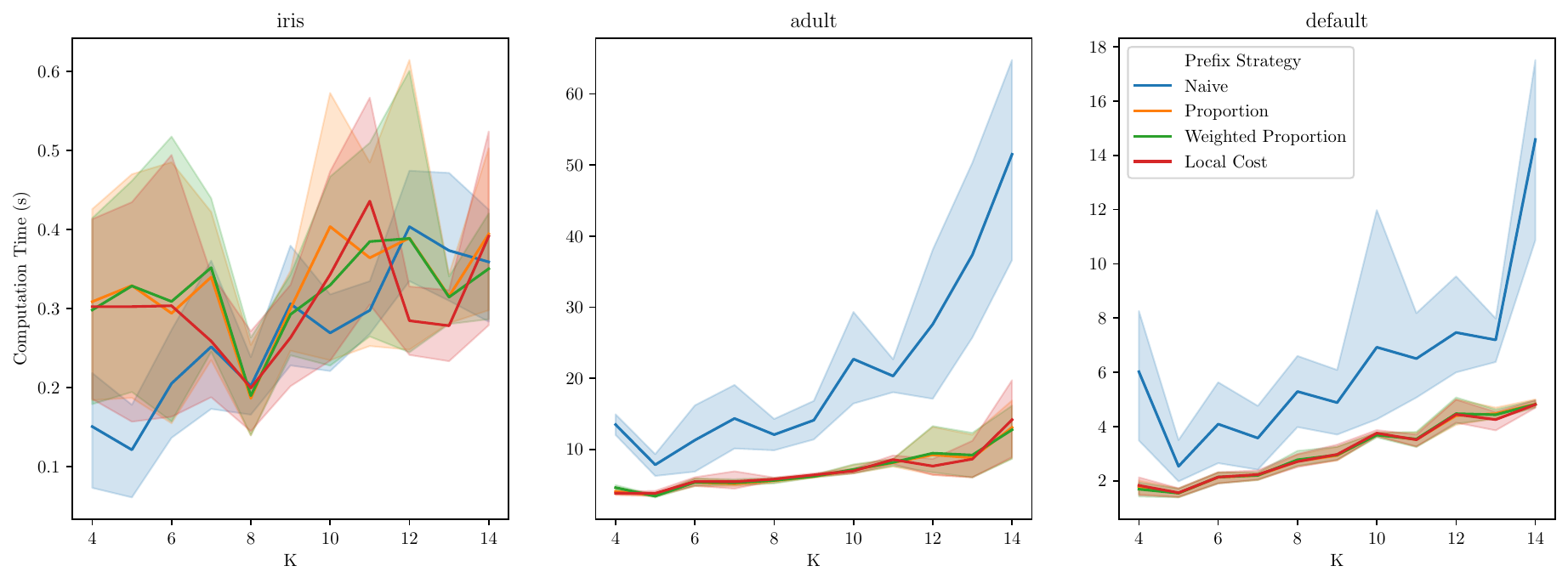}
    \end{subfigure}
  \caption {\label{fig:prefix_time} Impact of pre-fixing strategy on the run-time of the MiniReL algorithm. Proportion, Weighted Proportion, and Local Cost are pre-fixing strategies that use the pre-fix IP model with different objectives. Naive denotes random assignments of groups to clusters. Results averaged over 10 random seeds, widths indicate standard error.
  }
\end{figure}

Figure \ref{fig:prefix_iter} shows the impact of pre-fixing strategy on the number of iterations needed for MiniReL to converge. For iris, pre-fixing has practically no impact on on the number of iterations, however for larger datasets like adult and default using the pre-fix IP model with any objective leads to substantially fewer iterations. The same holds for cluster cost as shown in Figure \ref{fig:prefix_cost} where the IP model leads to solutions with slightly better clustering cost in both adult and default. Both results show that the choice of objective function has relatively little impact on the performance of pre-fixing, but outperform random assignment. 

\begin{figure}[!htb]
    \centering
    \begin{subfigure}
      \centering
      \includegraphics[width=\textwidth]{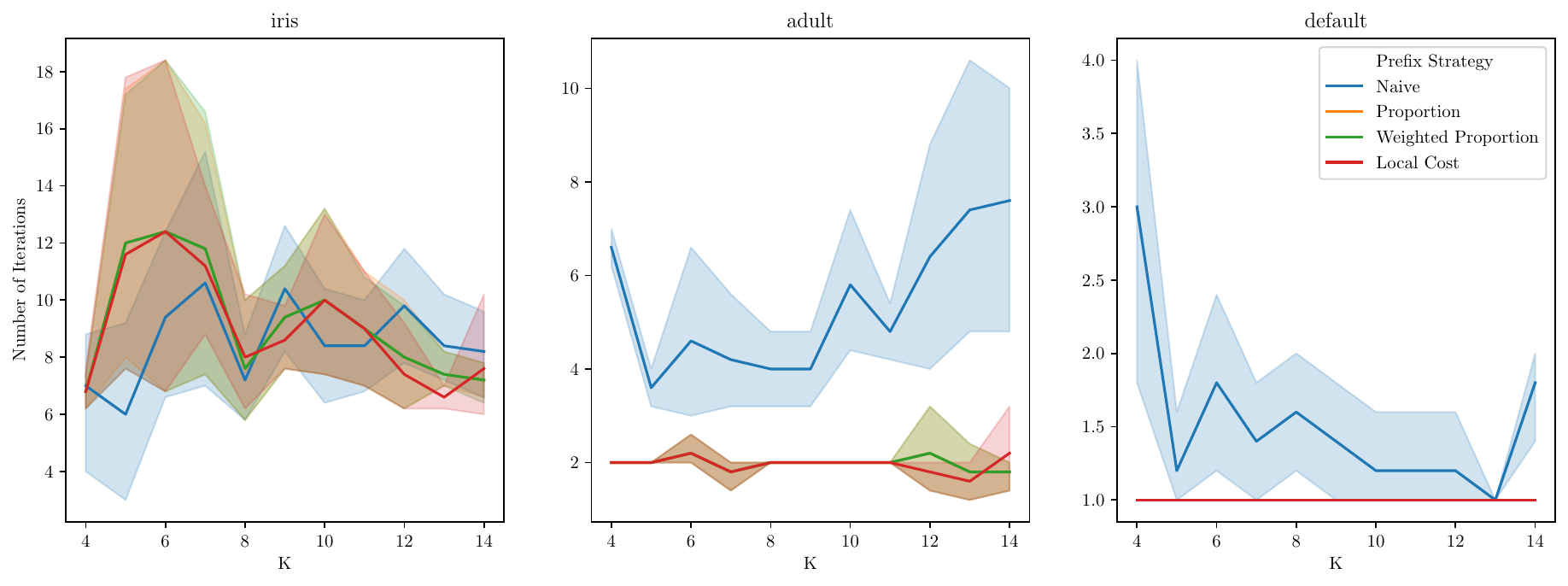}
    \end{subfigure}
  \caption {\label{fig:prefix_iter} Impact of pre-fixing strategy on the on the number of iterations to convergence for MiniReL. Results averaged over 10 random seeds. Bars indicate standard error.
  }
\end{figure}
\begin{figure}[!htb]
    \centering
    \begin{subfigure}
      \centering
      \includegraphics[width=\textwidth]{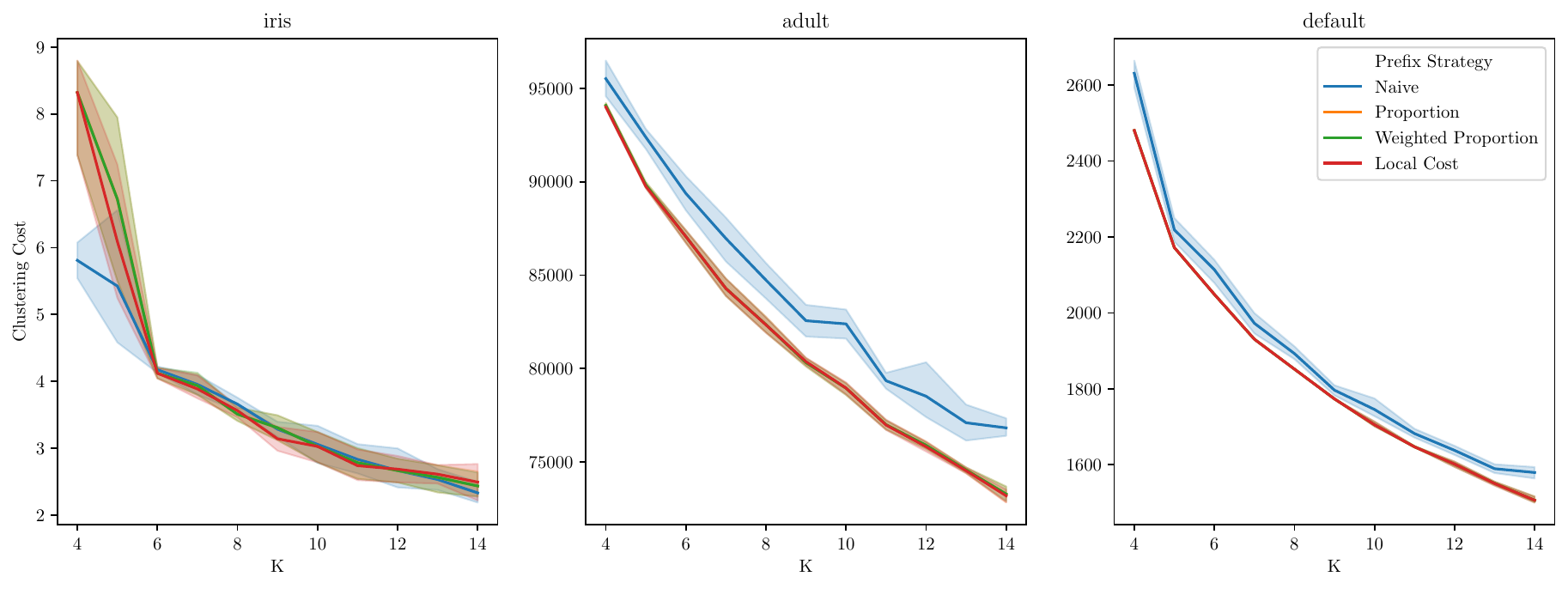}
    \end{subfigure}
  \caption {\label{fig:prefix_cost} Impact of pre-fix strategy on the clustering cost of the solution output by MiniReL. Results averaged over 10 random seeds. Bars indicate standard error.
  }
\end{figure}

\section{Warm-Starting MiniReL} \label{app:warm_start}
To reduce the number of iterations needed to converge in MiniReL, we warm-start the initial cluster centers with the final centers of the unfair variants of Lloyd's algorithm. The key intuition behind this approach is that it allows us to leverage the polynomial time assignment problem for the majority of iterations, and only requires solving the fair assignment problem to adjust the locally optimal unfair solution to a fair one. To incorporate warm-starting into MiniReL, we replace step 1 in Algorithm \eqref{alg:minrepfairlloyd} with centers generated from running Lloyd's algorithm. For the sake of brevity, we report results in the $k$-means setting. We benchmark this approach against two baselines: $(i)$ randomly sampling the center points, and $(ii)$ using the $k$-means++ initialization scheme without running Lloyd's algorithm afterwards. We also compare using the $k$-means warm-start with 1 initialization and 100 initializations. Figure \ref{fig:init_time} shows the impact of these initialization schemes on the total computation time including time to perform the initialization. Each initialization scheme was tested on three datasets from the UCI machine learning repository \cite{dua2017uci}: Iris, Adult, and Default. For each dataset we  randomly sub-sample 2000 data points (if $n > 2000$), and re-run MiniReL with 10 random seeds. The results show that using Lloyd's algorithm to warm-start MiniReL can lead to a large reduction in computation time, even taking into account the cost of running the initialization. However, there are diminishing returns. Namely running 100 different initialization for $k$-means and selecting the best leads to slightly larger overall run-times. 

\begin{figure}[t]
    \centering
    \begin{subfigure}
      \centering
      \includegraphics[width=.8\textwidth]{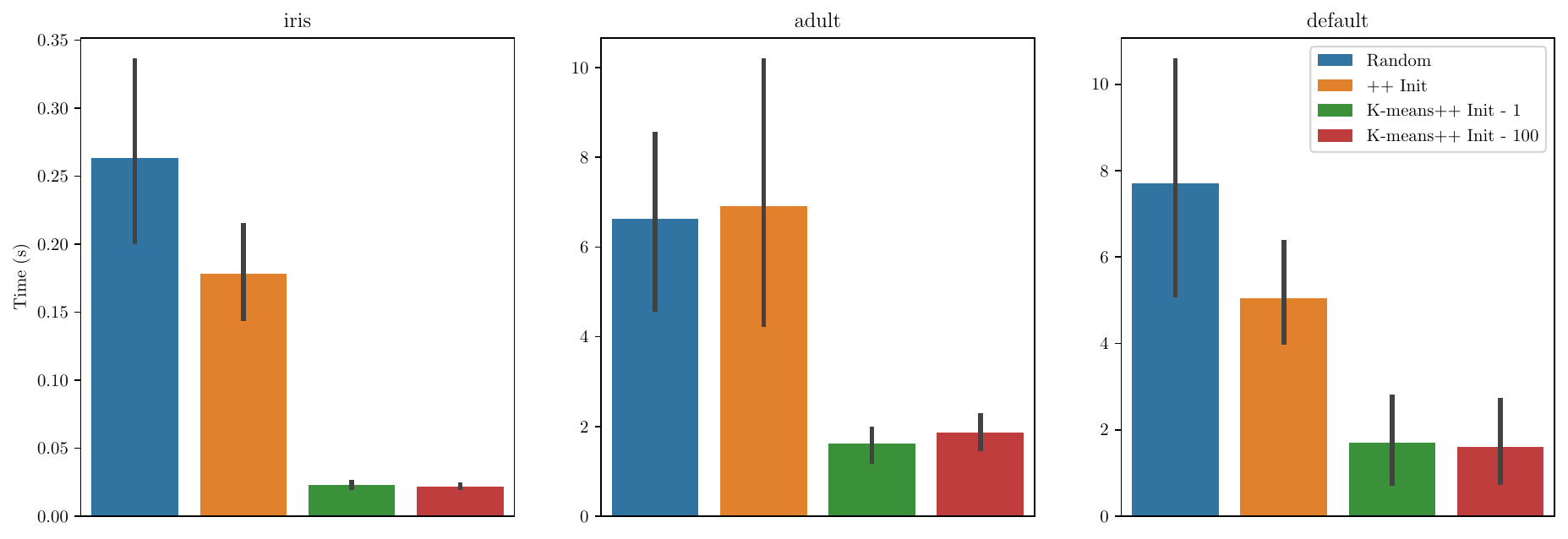}
    \end{subfigure}
  \caption {\label{fig:init_time} Average computation time in seconds over 10 random seeds for MiniReL with different initialization schemes. Bars indicate standard error.
  }
\end{figure}		
\section{Experiments with different values of $\alpha$} 
\label{app:changing_alpha}
In this section we explore the impact of $\alpha$ on the fairness and runtime of the different algorithms. For these experiments we replicate the experimental set-up of Section \ref{sec:kmeans_exp} but vary $\alpha$. For brevity we focus on the $k$-means setting, but $k$-medians showed similar results. Figures \ref{fig:alpha_fv} and \ref{fig:alpha_add_fv} show the normalized fairness violation (the number of $\alpha$-represented clusters away from meeting the required amount $\beta_g$ normalized by $K$) and normalized additive fairness violation (the number of data points away from meeting the fairness constraints normalized by $n$). In both cases, lower value of $\alpha$ lead to much laxer fairness requirements that are easier to meet. Both the standard $k$-means and the heuristic fair $k$-means algorithms are able to generate fair clusters for $\alpha \leq 0.3$ for both statistical parity and equality of opportunity settings of $\beta_g$. However for larger values of $\alpha$ both struggle to generate fair clusters having normalized (additive) violations close to 100\% (25\%). In constrast, by design MiniReL with IP-based assignment (MiniRel-Prefix) is able to generate fair clusters at all levels of $\alpha$. While the version of MiniReL with network flow rounding (MiniReL-PrefixFlow) can lead to fairness violations, the size of the violations are consistently small (in line with our theoretical results) - violating the constraints by close to 0\% of the size of the dataset. Figure \ref{fig:alpha_time} shows the impact of $\alpha$ and $K$ on the computation time of all algorithms. As expected, both $k$-means and fair $k$-means run very quickly regardless of $K$ (and $\alpha$ evidently has no impact on the runtime). In contrast, the runtime of MiniReL increases with both $K$ and $\alpha$ showing that harder assignment problems lead to larger runtimes overall.

\begin{figure}[!htb]
    \centering
    \begin{subfigure}
      \centering
      \includegraphics[width=\textwidth]{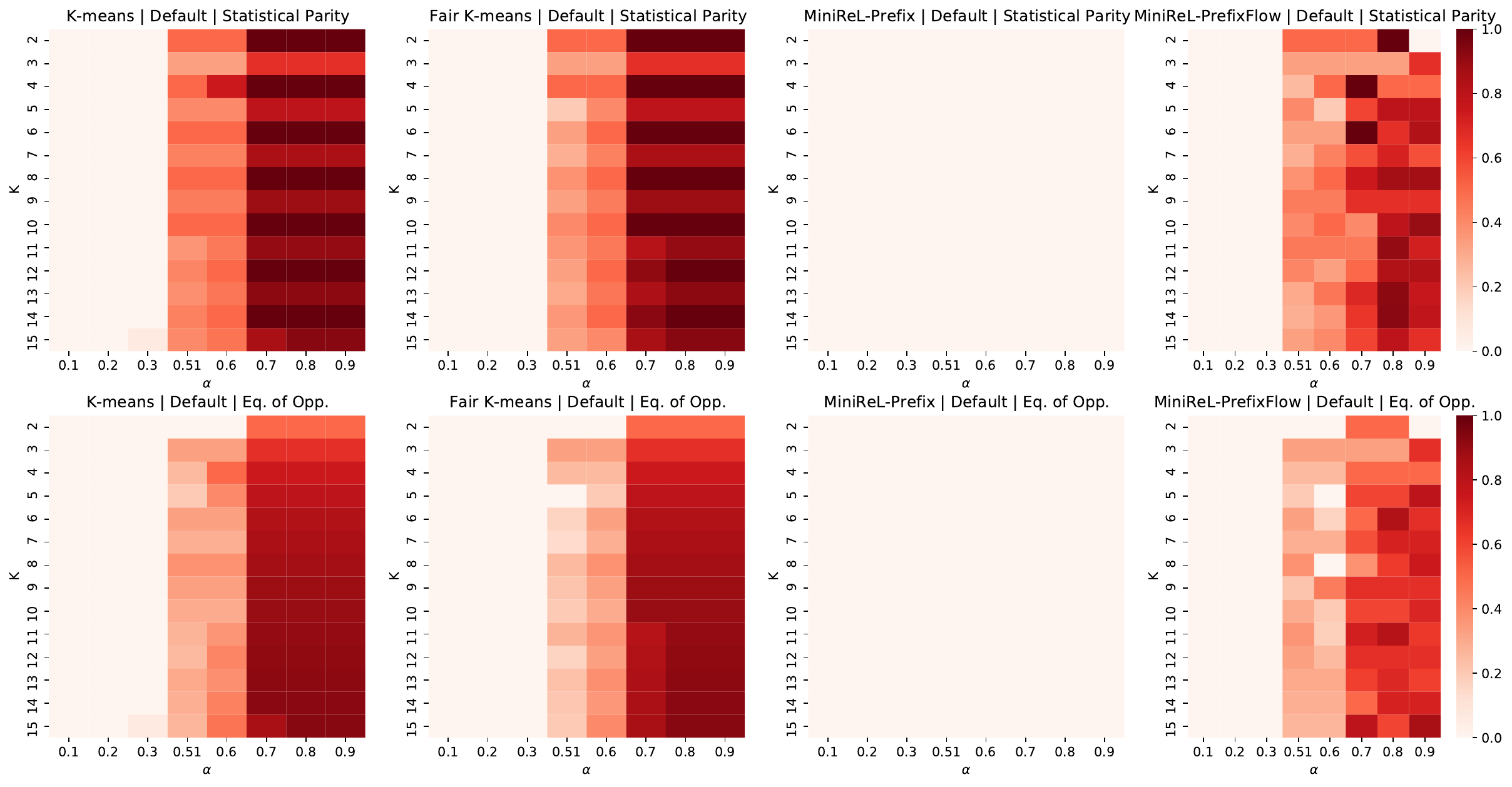}
    \end{subfigure}
  \caption {\label{fig:alpha_fv} Impact of $\alpha$ and $K$ on the normalized fairness violation ($\max_g\max(\beta_g - \Lambda(C, X_g, \alpha),0)/K$). Higher values (closer to one) indicate larger fairness violations.
  }
\end{figure}

\begin{figure}[!htb]
    \centering
    \begin{subfigure}
      \centering
      \includegraphics[width=\textwidth]{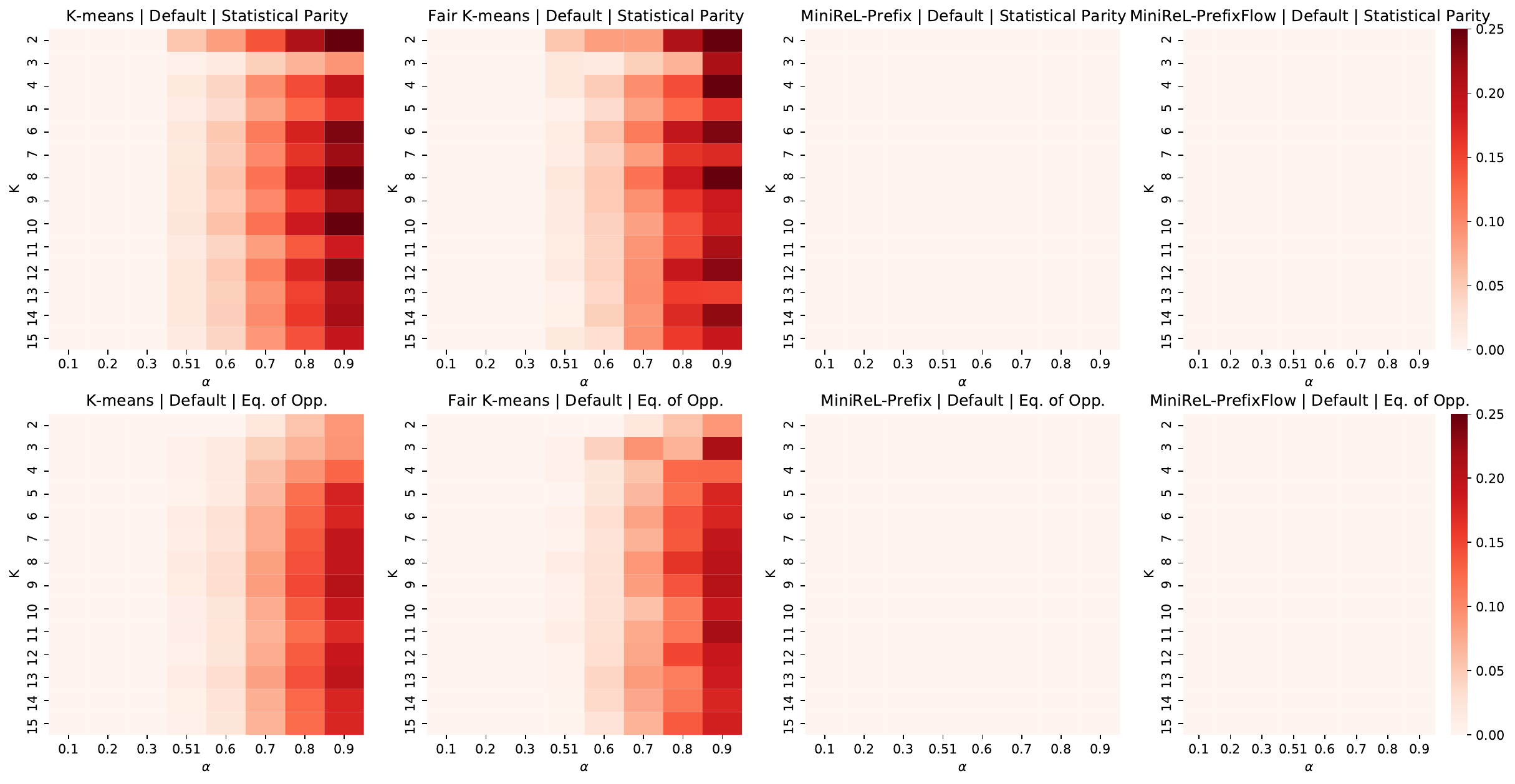}
    \end{subfigure}
  \caption {\label{fig:alpha_add_fv} Impact of $\alpha$ and $K$ on the normalized additive fairness violation (normalized by number of data points). Higher values (closer to one) indicate larger fairness violations.
  }
\end{figure}

\begin{figure}[!htb]
    \centering
    \begin{subfigure}
      \centering
      \includegraphics[width=\textwidth]{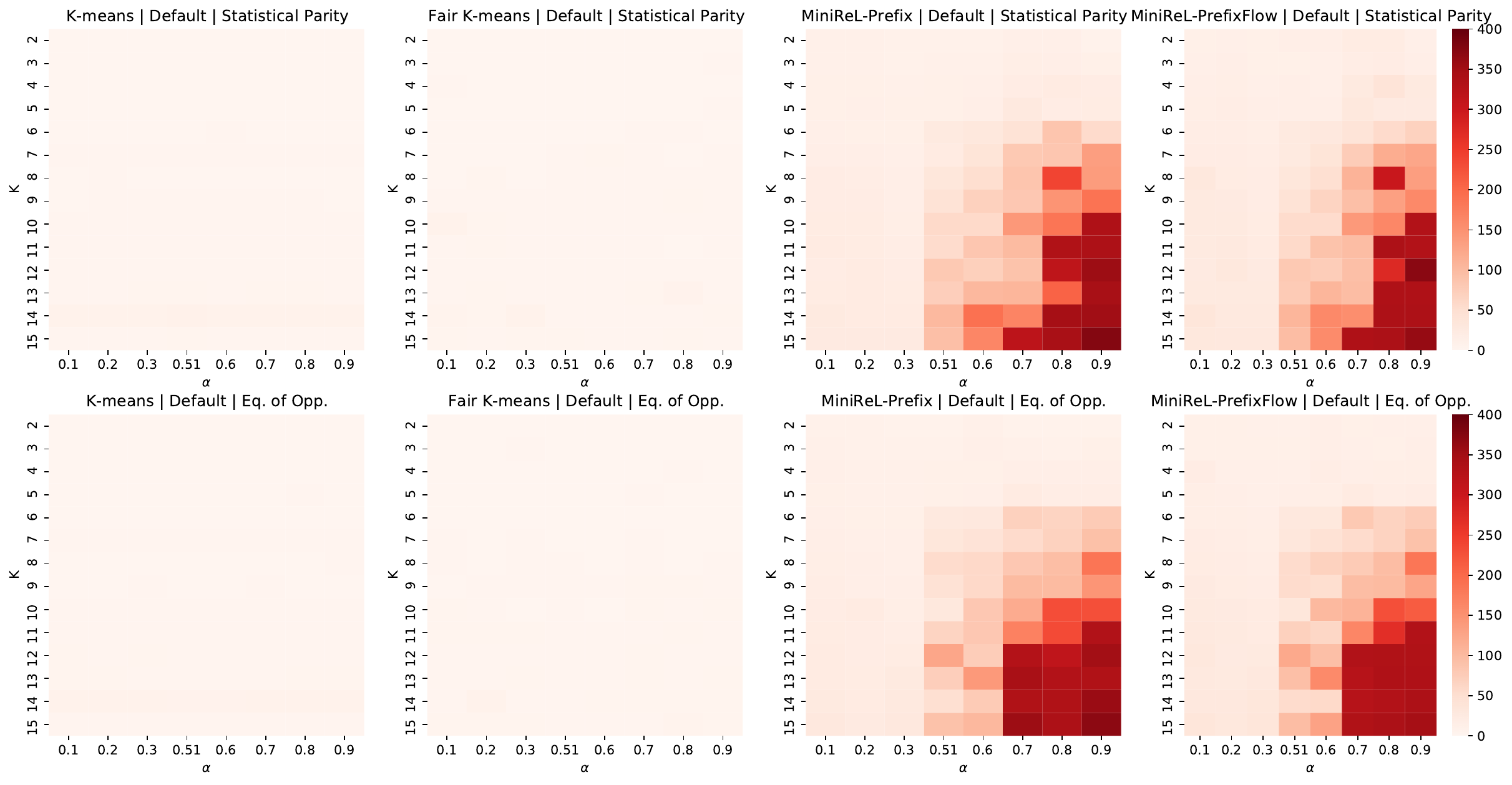}
    \end{subfigure}
  \caption {\label{fig:alpha_time} Impact of $\alpha$ and $k$ on the runtime (s) of the different algorithms.
  }
\end{figure}
		
\section{Experiments with balanced clusters} 
\label{app:balanced_clusters}
In this section we explore the impact of constructing \textit{balanced clusters} (i.e., clusters with approximately equal size) on the objective and computation time of MiniReL. For these experiments we require each cluster to have at least 80\% of the balanced allocation of data points to clusters (i.e., $\ell = 0.8*\frac{n}{K}$. We compare their performance on three large datasets: Adult, Default, and Voting. For the Adult dataset we include one instance where we use one sensitive feature (gender), and one where we use two (gender and marital status). Figure \ref{fig:balanaced_obj} shows the $k$-means clustering cost of both Lloyd's algorithm (with no cardinality constraints) and the cost of MiniReL with balanced clusters. The first thing to note is that under balanced clusters, not every dataset and setting for $\alpha$, $\beta$ remains feasible. For the voting dataset under statistical parity, MiniReL certifies that no fair clustering exists (and thus has no associated objective or runtime). Overall, requiring balance leads to a modest increase in cluster cost across all instances - exceeding the deviation that occurs from enforcing fairness alone. Table \ref{tab:balanced_sp_runtime} shows the average runtime over 10 random seeds for $K$-Means and MiniReL-HeurFlow under balanced clusters for statistical parity. Adding the requirement that clusters are balanced leads to a small increase in runtime for MiniReL on all datasets. The one exception is the voting dataset, which is an infeasible problem, where MiniReL can certify infeasibility in under 30 seconds for all settings of $K$.

\begin{figure}[!htb]
    \centering
    \begin{subfigure}
      \centering
      \includegraphics[width=\textwidth]{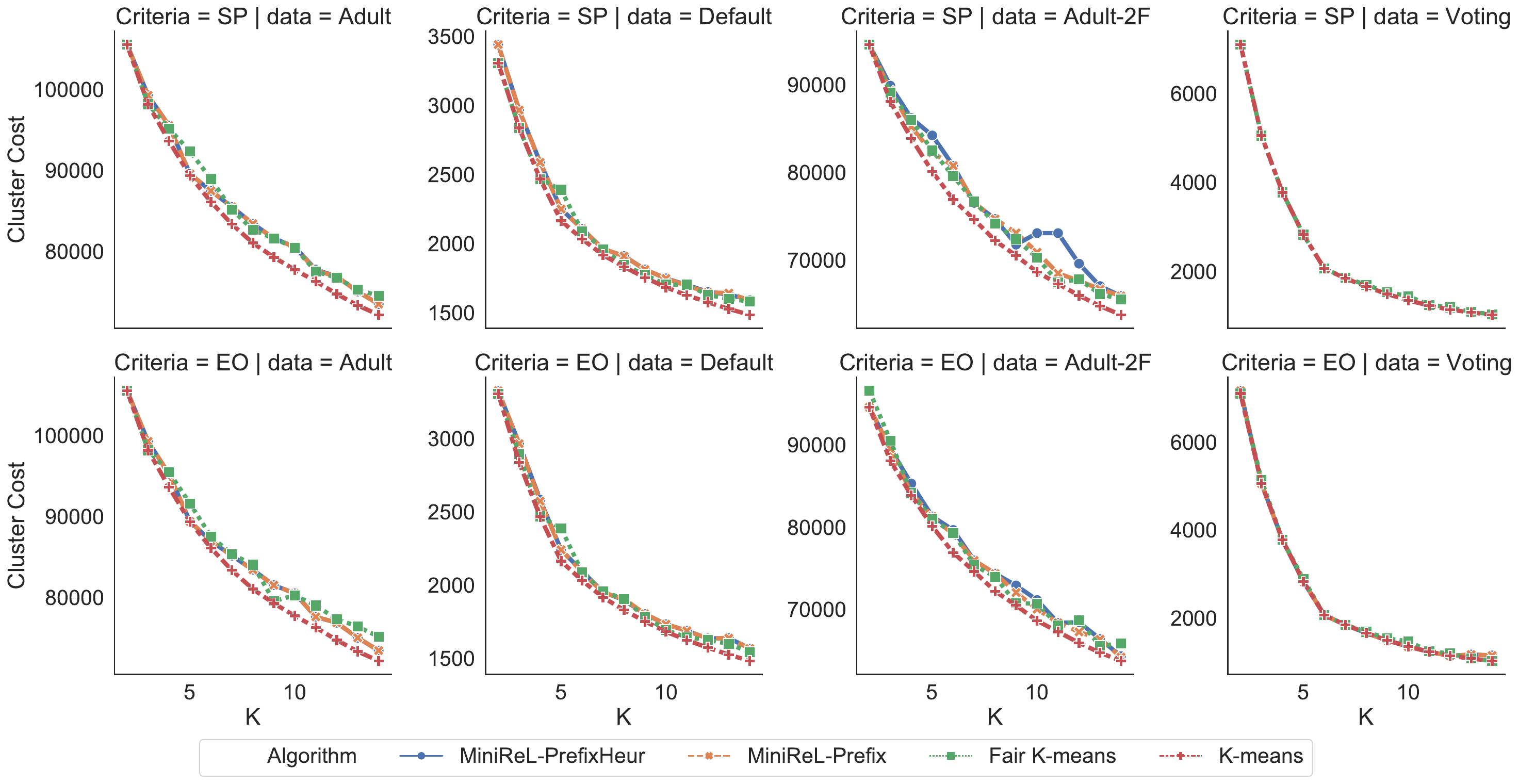}
    \end{subfigure}
  \caption {\label{fig:balanaced_obj} $k$-means clustering cost of Lloyd's algorithm and MiniReL with balanced clusters.
  }
\end{figure}

\begin{table}
\footnotesize
\centering
\caption{Average computation time in seconds for MiniReL-HeurFlow under Statistical Parity with balanced clusters. Final row includes average computation time over all settings of K for each dataset. \label{tab:balanced_sp_runtime}}.
\begin{tabular}{lrrrrrrrr}
\toprule
data & \multicolumn{2}{c}{Adult} & \multicolumn{2}{c}{Adult-2F} & \multicolumn{2}{c}{Default} & \multicolumn{2}{c}{Voting} \\
k & K-means & MiniReL &  K-means & MiniReL & K-means & MiniReL & K-means & MiniReL \\
\midrule
2    &    2.22 &             9.41 &     3.04 &             8.80 &    0.65 &            11.21 &    0.48 &             8.64 \\
3    &    2.40 &            15.08 &     3.03 &            13.47 &    0.82 &            11.72 &    0.56 &             9.23 \\
4    &    2.80 &            26.29 &     3.00 &            27.71 &    0.79 &            12.79 &    1.30 &            10.24 \\
5    &    2.58 &            18.27 &     3.59 &            30.41 &    1.25 &            13.85 &    0.84 &            13.31 \\
6    &    2.87 &            19.00 &     2.65 &            22.79 &    1.42 &            15.03 &    0.92 &            13.22 \\
7    &    3.10 &            39.18 &     3.66 &            23.63 &    2.36 &            16.10 &    1.07 &            14.37 \\
8    &    3.08 &            41.11 &     5.70 &            25.51 &    1.50 &            17.83 &    1.48 &            15.27 \\
9    &    3.15 &            29.61 &     6.07 &            76.04 &    1.53 &            18.57 &    1.53 &            13.64 \\
10   &    4.57 &            31.00 &     6.33 &            70.21 &    1.88 &            20.45 &    2.76 &            15.12 \\
11   &    4.02 &            24.87 &     3.79 &            46.01 &    2.01 &            33.05 &    1.90 &            15.94 \\
12   &    3.53 &            42.30 &     5.24 &            72.27 &    2.66 &            34.49 &    2.93 &            17.71 \\
13   &    4.09 &            36.90 &     4.08 &            81.45 &    3.21 &            30.42 &    3.23 &            23.77 \\
14   &    4.33 &            27.62 &     2.91 &            46.12 &    8.42 &            24.61 &    3.21 &            22.37 \\ \midrule
mean &    3.29 &            27.74 &     4.08 &            41.88 &    2.19 &            20.01 &    1.71 &            14.83 \\
\bottomrule
\end{tabular}

\end{table}

\section{Full $k$-means results} \label{app:kmean_extra}
The following section includes the extended empirical results for MiniReL in the $k$-means setting. Figure \ref{fig:kmeans_eo} shows both the maximum fairness violation and normalized additive fairness violation of both algorithms under equality of opportunity. 

\begin{figure}[t]
    \centering
    \begin{subfigure}
      \centering
      \includegraphics[width=0.9\textwidth]{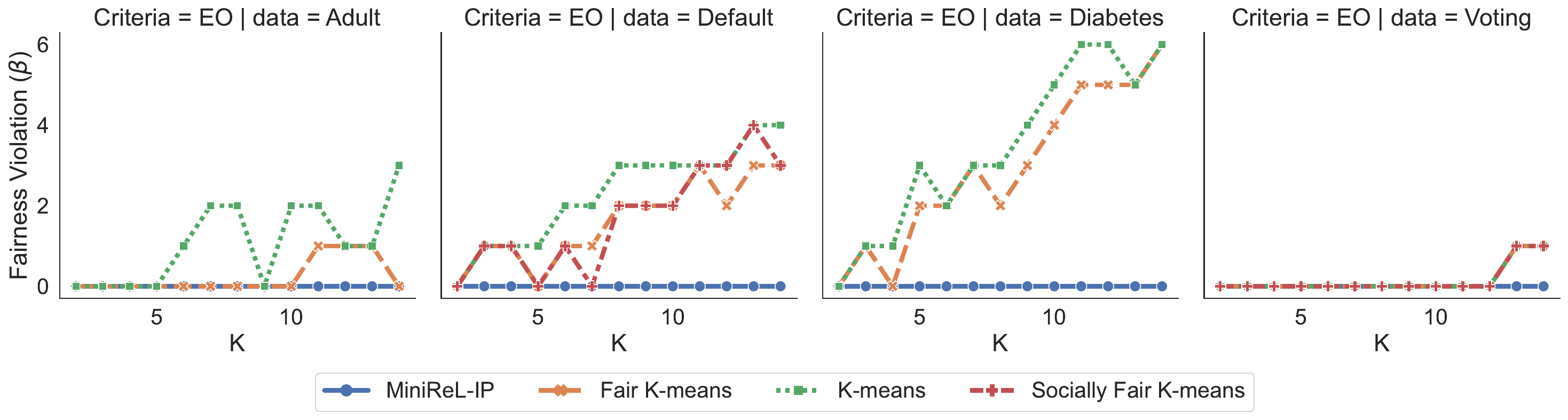}
    \end{subfigure}
    \begin{subfigure}
      \centering
      \includegraphics[width=0.9\textwidth]{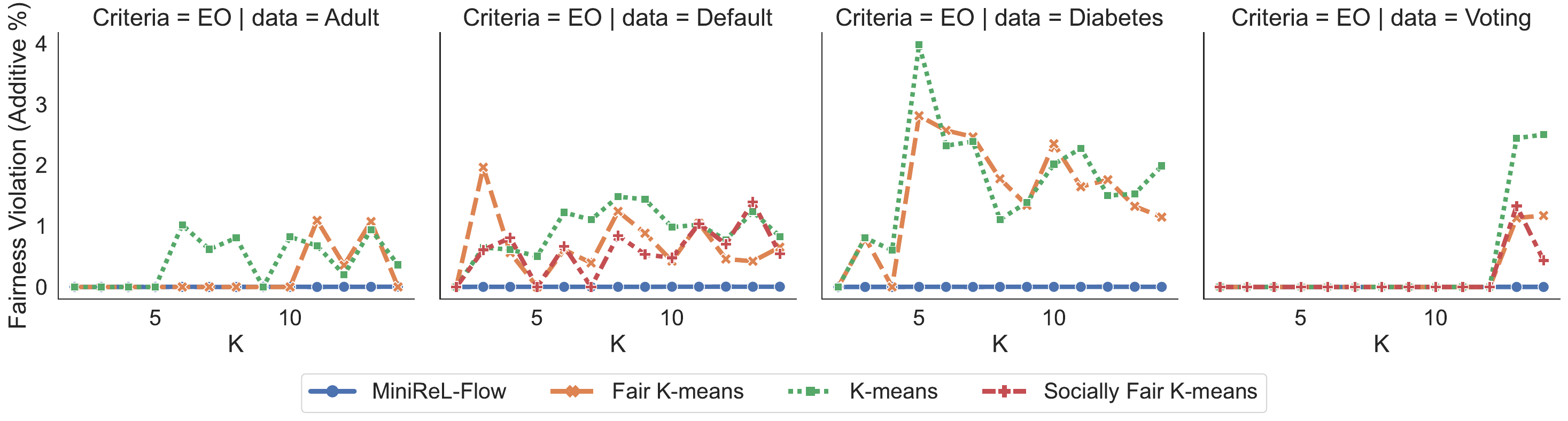}
    \end{subfigure}
  \caption {\label{fig:kmeans_eo} \vspace{-0.1cm}(Top) Fairness violation ($\beta_g$) for MiniReL-IP and baseline algorithms in $k$-means setting under equality of opportunity. (Bottom) Normalized Additive Fairness Violation (additive fairness violation divided by dataset size) for MiniReL and baseline algorithms under equality of opportunity.
  }
\end{figure}

\begin{figure}[t]
    \centering
    \begin{subfigure}
      \centering
      \includegraphics[width=\textwidth]{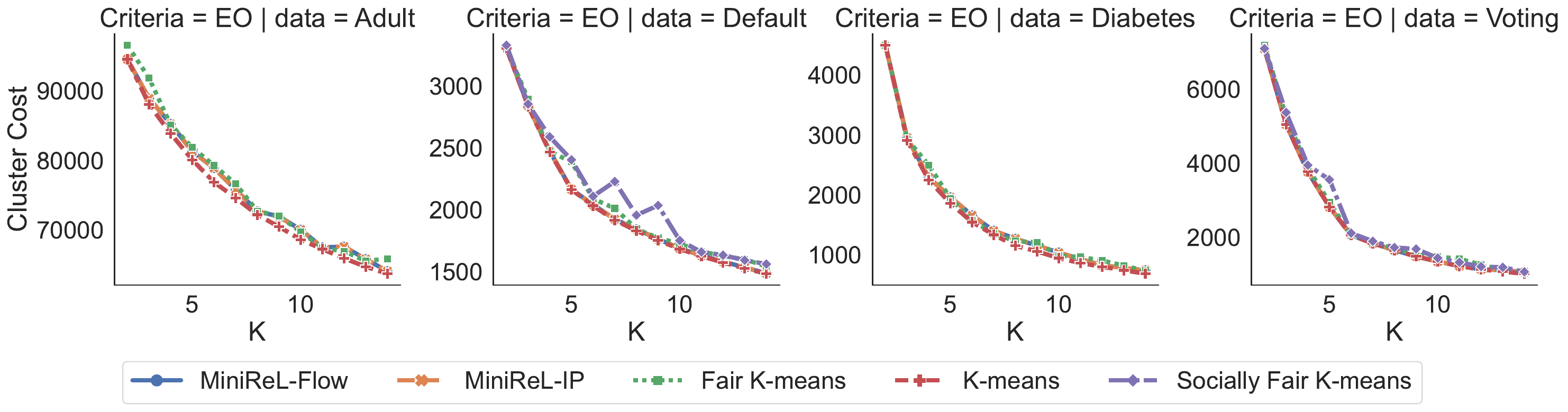}
    \end{subfigure}
    \vspace{-0.5cm}
  \caption {\label{fig:cluster_cost_kmeans_eo} $k$-means clustering cost of fair clustering algorithms under equality of opportunity.
  }
\end{figure}		
\section{Full $k$-medians results} \label{app:kmed}
The following section includes the extended empirical results for MiniReL in the $k$-medians section. Figure \ref{fig:kmed_fairness_joint} shows both the maximum fairness violation and normalized additive fairness violation of both algorithms under both notions of fairness. Regardless of the setting, the baseline $k$-medians algorithm and its fair heuristic analog are unable to generate fair clusters - whereas MiniReL can construct fair clusters or clusters with negligible additive fairness violations by design. Figure \ref{fig:kmed_cluster_cost_joint} shows the clustering objective of both sets of algorithms. With the exception of the voting dataset under statistical parity, MiniReL is able to generate fair clusters with practically no increase in the cluster cost - showing that in practical settings achieving minimum representation fairness comes at little cost over Lloyd's algorithm.

\begin{figure}[t]
    \centering
    \begin{subfigure}
      \centering
      \includegraphics[width=0.9\textwidth]{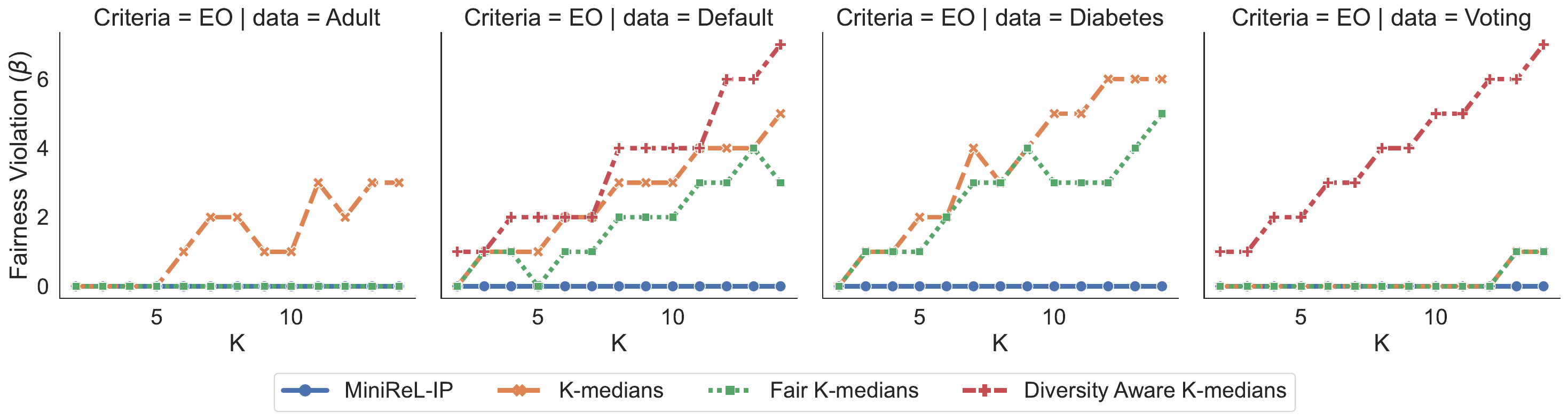}
    \end{subfigure}
    \begin{subfigure}
      \centering
      \includegraphics[width=0.9\textwidth]{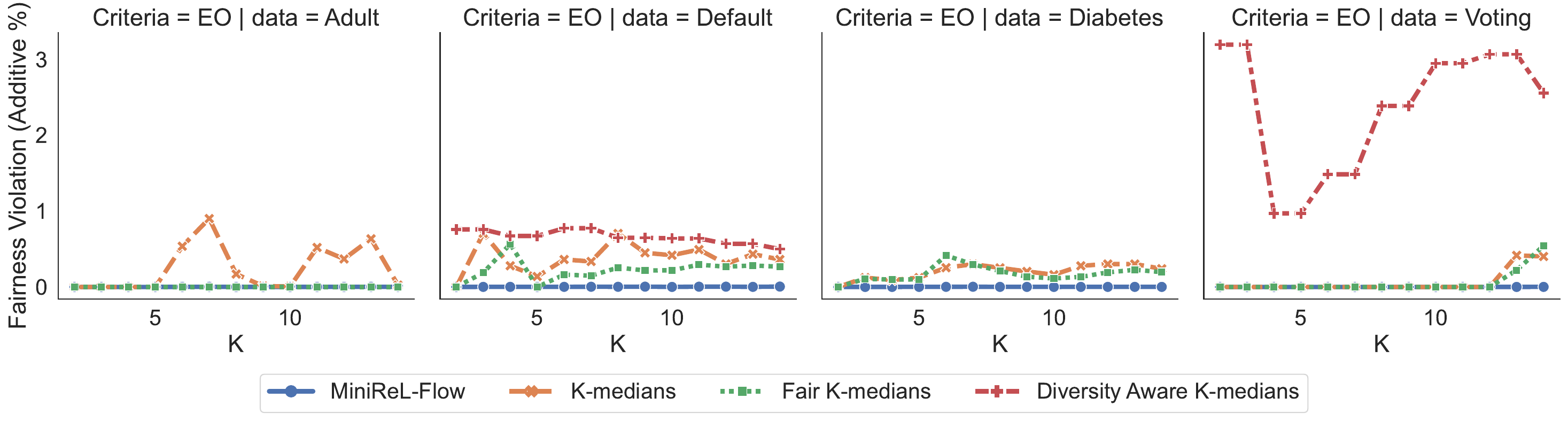}
    \end{subfigure}
  \caption {\label{fig:kmed_fairness_joint} \vspace{-0.1cm}(Top 2 Rows) Fairness violation ($\beta_g$) for baseline algorithms and MiniReL in the $k$-medians setting under equality of opportunity. (Bottom 2 Rows) Normalized Additive Fairness Violation (additive fairness violation divided by dataset size) for baseline algorithms and MiniReL under equality of opportunity.
  }
\end{figure}

\begin{figure}[t]
    \centering
    \begin{subfigure}
      \centering
      \includegraphics[width=\textwidth]{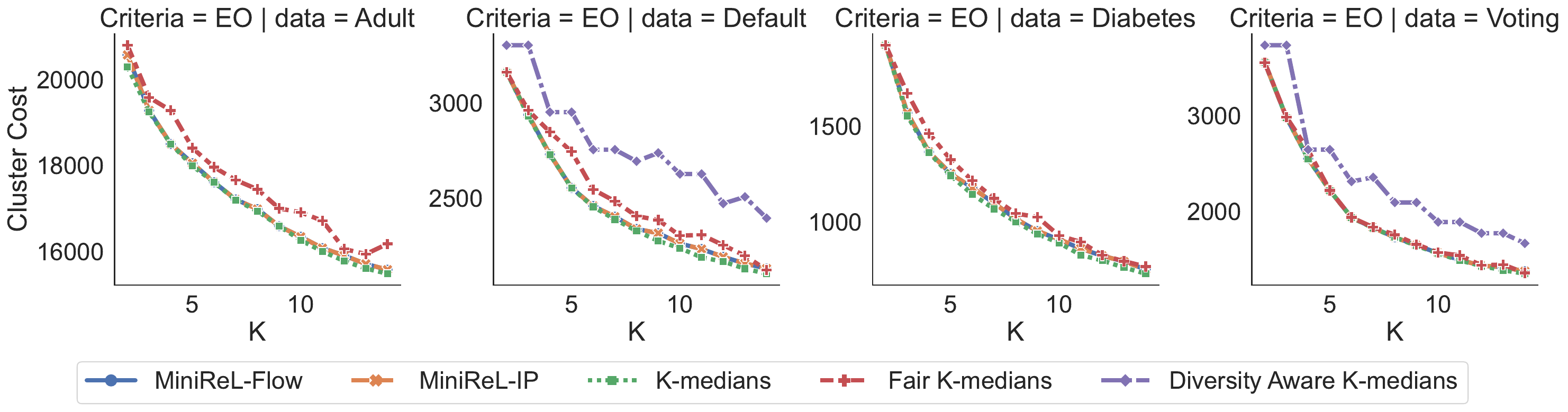}
    \end{subfigure}
    \vspace{-0.5cm}
  \caption {\label{fig:kmed_cluster_cost_joint} $k$-medians clustering cost of baseline algorithms and MiniReL under equality of opportunity.
  }
\end{figure}		

\section{Experiments with Synthetic Data} \label{app:synthetic}

To investigate how MiniReL scales with respect to various problem dimensions we benchmark the algorithm on synthetic data. We generate the dataset using an isotropic gaussian mixture model included in scikit-learn \cite{scikit-learn}. We vary the dimension of the data ($d$) and the number of data points ($n$). For all experiments we fix $k=10$. To generate sensitive features and labels we sample one of ${|\cal G}|$ labels uniformly at random, and run this labeling process independently $|\cal F|$ times (once for each feature). For each experiment we fix all but one of $n=10000$, $d=100$, $|{\cal F}=1|$ and $|{\cal G}|=2$ and evaluate the impact of changing the remain dimension. \ref{fig:synthetic} summarizes the results for three variants of MiniReL: MiniReL-HeurFlow, MiniReL-Prefix, and MiniReL-PrefixFlow. 

The results indicate that the dimension of the underlying data $d$ has a negligible impact on the runtime of MiniReL. This is to be expected as the dimension of the data only impacts the computation of the objective in the FMRA which is a small fraction of the overall computation time. The runtime of MiniReL does scale with the size of the dataset as $n$. Notably, the heuristic pre-fixing algorithm outperforms IP-based methods for solving the RAP as the dataset size grows. Likewise, it scales better with respect to both the number of sensitive features and the number of groups.

\begin{figure}[!htb]
    \centering
    \begin{subfigure}
      \centering
      \includegraphics[width=0.45\textwidth]{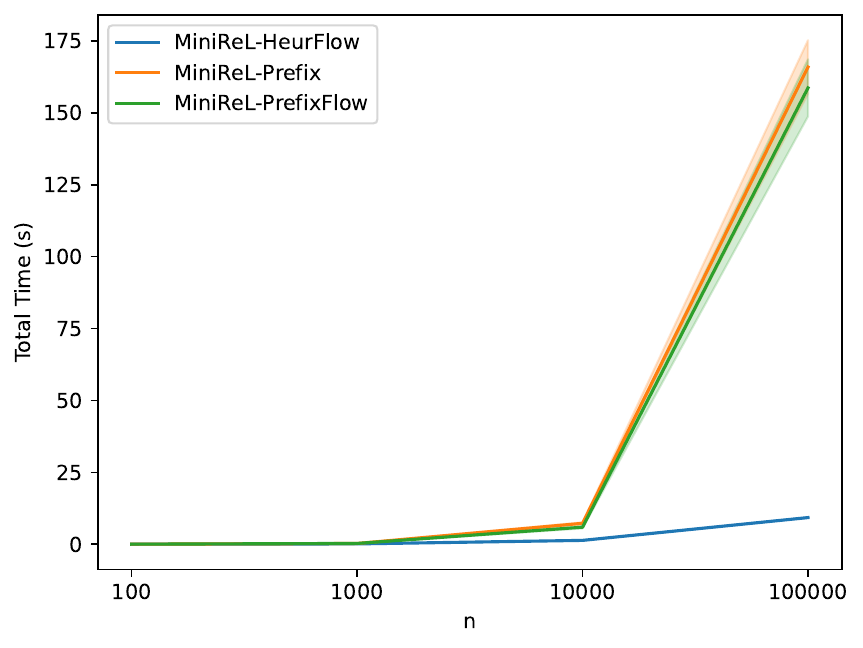}
    \end{subfigure}
        \begin{subfigure}
      \centering
      \includegraphics[width=0.45\textwidth]{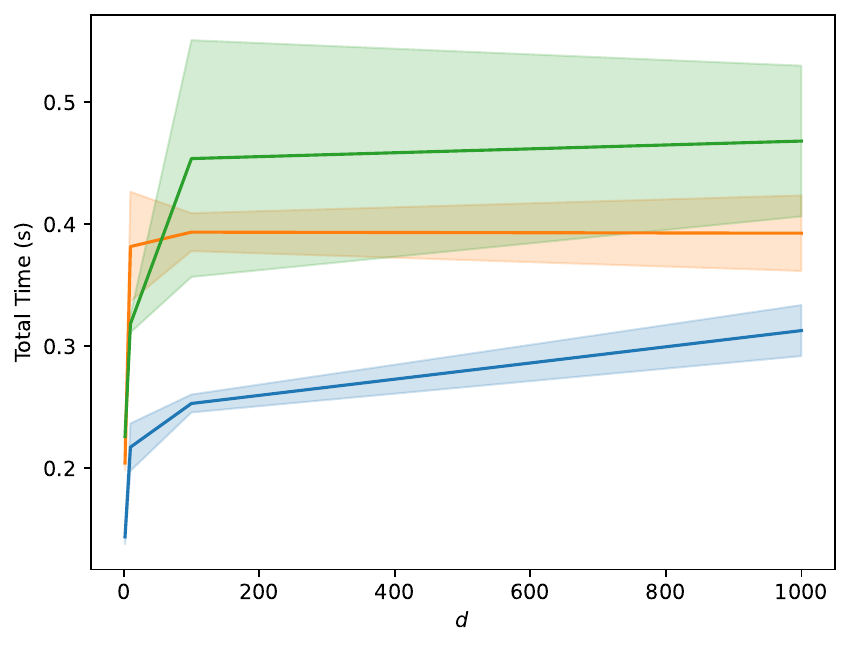}
    \end{subfigure}
    \begin{subfigure}
      \centering
      \includegraphics[width=0.45\textwidth]{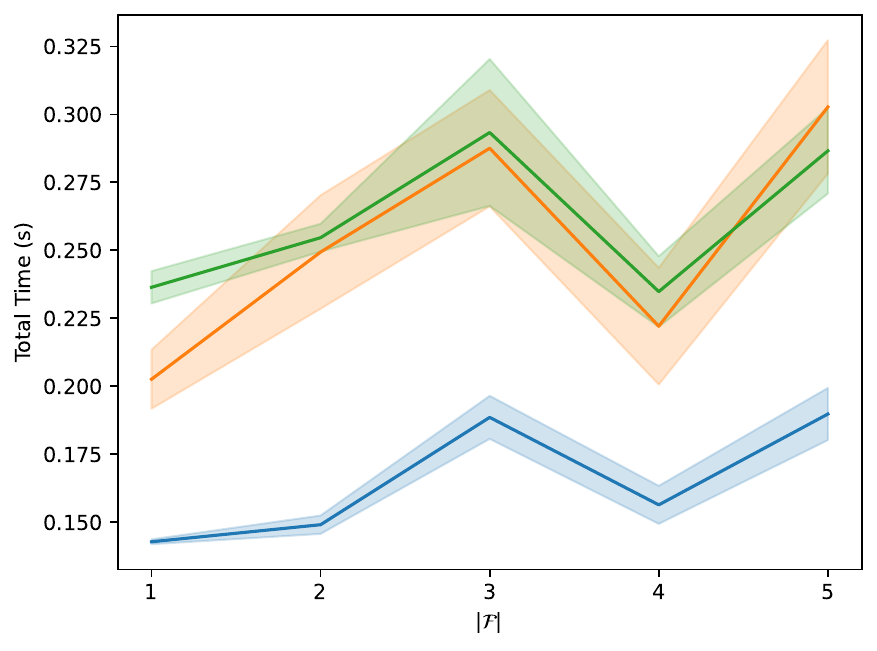}
    \end{subfigure}
    \begin{subfigure}
      \centering
      \includegraphics[width=0.45\textwidth]{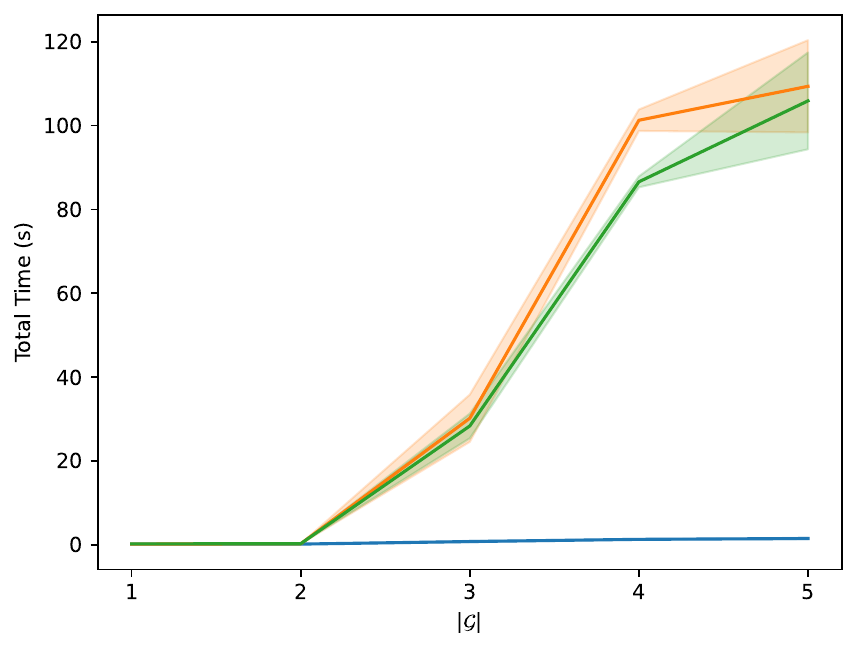}
    \end{subfigure}

  \caption {\label{fig:synthetic} Impact of data set size ($n$), number of features ($d$), number of sensitive features ($|{\cal F}|$), and number of groups ($|{\cal G}|$) on the runtime of MiniReL algorithms on synthetic data.
  }
\end{figure}

\end{APPENDIX}

	\end{document}